\documentclass[a4paper,10pt]{article}
\usepackage{amsmath}
\usepackage{amsfonts}
\usepackage{amssymb}
\usepackage{latexsym}
\usepackage{epsfig}
\usepackage{graphicx}
\usepackage{oldgerm}
\usepackage{theorem}

\def\0{{\bf 0}}

\def\g{g^{(\nu,\kappa)}_{N,T}}
\def\p{\widetilde{p}^{(\nu,\kappa)}}
\def\pT{p_T^{(\nu,\kappa)}}

\def\x{{\bf x}}
\def\y{{\bf y}}

\def\ma{{\mathfrak a}}
\def\mb{{\mathfrak b}}

\def\mg{\mathfrak{g}^{(\nu,\kappa)}_{N,T}}
\def\mp{\mathfrak{p}^{(\nu,\kappa)}_{N,T}}

\def\Knu{\widetilde{J}_\nu}
\def\hL{\widehat{L}}
\def\hKnu{\widehat{J}_\nu}

\def\N{\mathbb{N}}
\def\R{\mathbb{R}}
\def\Rp{\mathbb{R}_+}
\def\mX{\mathfrak{X}}

\def\X{{\bf X}}
\def\Y{{\bf Y}}
\def\Z{\mathbb{Z}}

\def\tS{\tilde{S}}
\def\tI{\tilde{I}}

\def\hR{\widehat{R}}
\def\hPhi{\widehat{\Phi}}

\def\CT{C_{N,T}^{\nu,\kappa}}

\def\RW{\mathbb{R}^{N}_{+<}}

\def\bA{\mathbb{A}}
\def\bS{\mathbb{S}}

\def\ba{{\bf a}}

\def\cD{{\cal D}}
\def\cI{\tilde{\cal I}}
\def\cS{{\cal S}}
\def\tcS{\widetilde{\cal S}}
\def\cG{{\cal G}}
\def\cR{{\rho}^{\bf Y}}
\def\cP{{\rho}^{\bf X}}

\def\Pf{{\rm Pf}}

\def\valpha{\mib{\alpha}}
\def\vbeta{\mib{\beta}}

\theorembodyfont{\itshape}

\newtheorem{thm}{Theorem}[section]
\newtheorem{lem}[thm]{Lemma}

\newtheorem{prop}[thm]{Proposition}

\newcommand{\mib}[1]{\mbox{\boldmath $#1$}}
\newcommand{\SSC}[1]{\section{#1}\setcounter{equation}{0}}
\newcommand{\qed}{\hbox{\rule[-2pt]{3pt}{6pt}}}

\def\vtheta{\mib{\theta}}
\def\kPsi{\Psi^{\bf Y}_{N,T}}
\def\kZ{Z^{\bf Y}_{N,T}}
\def\E{\mathbb{E}}
\def\D0{_{0}{\bf D}}
\def\Dinf{{\bf D}_{\infty}}

\begin{document}
\pagestyle{plain}

\noindent
{\bf \Large{Infinite Systems of 
Non-Colliding Generalized Meanders \\
and Riemann-Liouville Differintegrals} }

\vskip 0.5cm

\vskip 0.5cm

\noindent 
{\bf \large{Makoto Katori$^{1*}$, \, Hideki Tanemura$^{2**}$}} \\ 

\begin{description}
\item{1} \quad Department of Physics,
Faculty of Science and Engineering,
Chuo University, \\
Kasuga, Bunkyo-ku, Tokyo 112-8551, Japan. 
e-mail: katori@phys.chuo-u.ac.jp
\item{2} \quad
Department of Mathematics and Informatics,
Faculty of Science, Chiba University, \\
1-33 Yayoi-cho, Inage-ku, Chiba 263-8522, Japan.
e-mail:tanemura@math.s.chiba-u.ac.jp
\end{description}
\vskip 0.5cm


\vskip 0.5cm

\noindent
{\bf Abstract.} 
Yor's generalized meander is a 
temporally inhomogeneous modification of
the $2(\nu+1)$-dimensional Bessel process with $\nu > -1$,
in which the inhomogeneity is indexed by 
$\kappa \in [0, 2(\nu+1))$.
We introduce the non-colliding particle systems
of the generalized meanders and prove that 
they are Pfaffian processes, in the sense
that any multitime correlation function is
given by a Pfaffian. 
In the infinite particle limit, 
we show that the elements of matrix kernels
of the obtained infinite Pfaffian processes
are generally expressed by 
the Riemann-Liouville differintegrals 
of functions comprising the Bessel functions $J_{\nu}$
used in the fractional calculus,
where orders of differintegration 
are determined by $\nu-\kappa$.
As special cases of the two parameters
$(\nu, \kappa)$,
the present infinite systems include 
the quaternion determinantal processes 
studied by Forrester, Nagao and Honner
and by Nagao, which exhibit 
the temporal transitions
between the universality classes of
random matrix theory.
\vspace{0.5cm}

\noindent
{\bf running head:} Non-colliding generalized meanders

\vskip 0.3cm

\noindent
$^{*}$ Research supported in part by
the Grant-in-Aid for Scientific Research 
(KIBAN-C, No.17540363) of Japan Society for
the Promotion of Science.

\noindent
$^{**}$ Research supported in part by
the Grant-in-Aid for Scientific Research 
(KIBAN-C, No.15540106) of Japan Society for
the Promotion of Science.
\vskip 0.3cm

\vskip 0.3cm

\noindent
{\it Key words or phrases:} 
non-colliding generalized meanders, 
Bessel processes,
random matrix theory,
Fredholm Pfaffian and determinant,
Riemann-Liouville differintegrals,



\SSC{Introduction}\label{chap:Introduction}

The random matrix (RM) theory was introduced 
originally as an approximation theory of
{\it statistics} of nuclear energy levels \cite{Mehta04}.
It should be noted that at the same time
as the standard theory was established for three
ensembles called the Gaussian unitary, orthogonal,
and symplectic ensembles (GUE, GOE, GSE) \cite{Dys62b},
Dyson proposed to study such {\it stochastic processes}
of interacting particles that
the eigenvalue statistics of RMs are realized
in distribution of particle positions on $\R$
\cite{Dys62a}.
Dyson's Brownian motion model is a one-parameter
family of $N$-particle systems,
${\bf Z}^{(\beta)}(t)=(Z_{1}^{(\beta)}(t), 
Z_{2}^{(\beta)}(t), \cdots, Z_{N}^{(\beta)}(t))$, 
described by the stochastic differential equations
\begin{equation}
dZ_{i}^{(\beta)}(t)= dB_{i}(t)
+ \frac{\beta}{2} \sum_{1 \leq j \leq N, j \not=i}
\frac{1}{Z_{i}^{(\beta)}(t)-Z_{j}^{(\beta)}(t)} dt,
\quad t \in [0, \infty), 1 \leq i \leq N,
\label{eqn:Dyson1}
\end{equation}
where $B_{i}(t), i=1,2, \cdots, N$ 
are independent standard Brownian
motions and the parameter $\beta$ equals
2, 1 and 4 for GUE, GOE and GSE, respectively.
Due to the strong repulsive forces, which are
long-ranged and act between any pair of particles,
intersections of particle trajectories are
prohibited for $\beta \geq 1$ \cite{RS93}
(see also \cite{Bru89}).
In this one-parameter family, the $\beta=2$
case ({\it i.e.} the GUE case) is the simplest
and the most-understood, since its equivalence
with the $N$ particle systems of Brownian motions 
{\it conditioned never to collide with each other}
can be proved \cite{Gra99}.

The standard (Wigner-Dyson) theory has been extended 
by adding three chiral versions of RM ensembles
in the particle physics of QCD
\cite{VZ93,Ver94,JSV96,SV98}, 
and by introducing the four additional ensembles
so-called the Bogoliubov-de Gennes classes
in the mesoscopic physics \cite{AZ96,AZ97}.
Here we note that the chiral ensembles have
a parameter $\nu \in \{0,1,2, \cdots\}$
in addition to $\beta$.
In these totally ten ensembles \cite{AZ96,Z96,AZ97},
chiral GUE (chGUE), class C and class D
can be regarded as natural extensions of the GUE,
in the sense that these eigenvalue statistics
are also realized in appropriate non-colliding systems
of stochastic particle systems:
K\"onig and O'Connell showed that the chGUE with
the parameter $\nu \in \{0,1,2, \cdots\}$
corresponds to the non-colliding systems of
$2(\nu+1)$-dimensional squared Bessel processes
\cite{KO01}.
The present authors clarified that the 
eigenvalue statistics in the classes
C and D are realized by the non-colliding systems of
the Brownian motions {\it with an absorbing wall 
at the origin} and of 
the Brownian motions {\it reflecting at the origin}
\cite{KTNK03,KT04}.
Since the absorbing and reflecting Brownian motions
are directly related with the three-dimensional and
one-dimensional Bessel processes, respectively 
(see, for example, \cite{RY98}),
the stochastic differential equations of these
non-colliding particle systems are generally
given by
\begin{eqnarray}
&&d \widetilde{Z}_{i}^{(\nu)}(t)=dB_{i}(t)
+ \left[ \frac{2 \nu+1}{2} \frac{1}{\widetilde{Z}_{i}^{(\nu)}(t)}
+ \sum_{1 \leq j \leq N, j \not=i}
\left\{ \frac{1}{\widetilde{Z}_{i}^{(\nu)}(t)
-\widetilde{Z}_{j}^{(\nu)}(t)}
+ \frac{1}{\widetilde{Z}_{i}^{(\nu)}(t)
+\widetilde{Z}_{j}^{(\nu)}(t)}
\right\} \right] dt, \nonumber\\
&& \hskip 8cm 
t \in [0, \infty), 1 \leq i \leq N,
\label{eqn:Bessel}
\end{eqnarray}
with reflecting barrier condition at the origin in case $\nu=-1/2$.
Therefore, the difference of 
(non-standard) RM ensembles can be attributed to 
the difference of dimensionality of the Bessel processes,
whose non-colliding sets realize the statistics
of the RM ensembles \cite{KT04}.
Here we remind that the $d$-dimensional Bessel
process is defined as the process of the
radial coordinate (the modulus) of a Brownian motion
in $\R^{d}$.
To realize other $10-4=6$ RM ensembles
by conditioned stochastic processes
may be much more difficult
(see \cite{Spohn}), but we demonstrated that,
if we consider appropriate non-colliding systems
of {\it temporally inhomogeneous} processes
defined only in a finite time-interval $[0,T]$,
we can observe the transitions of distributions
into the 6 distributions as the time $t$ 
approaches the final time $T$ \cite{KT03b,KT04}.
The interesting fact is that the processes that 
can be used instead of the Bessel processes 
(\ref{eqn:Bessel}) should have one more parameter
$\kappa$ in addition to $\nu$.
This two-parameter family of temporally
inhomogeneous processes indexed by
$(\nu, \kappa),
\nu > -1, \kappa \in [0, 2(\nu+1))$
is equivalent with the family
of processes already studied by Yor.
He called them the {\it generalized meanders}
\cite{Yor92}.

From the view-point of random matrix theory,
studying time-development of stochastic systems
by calculating, for example, the multitime correlation
functions corresponds to considering
multi-matrix models. In particular,
the temporally inhomogeneous processes will be
identified with such matrix models
that matrices with different symmetries are 
coupled in a chain \cite{KT02a,KT03a,NKT03,KNT04}.
Determination of all multitime correlation functions
of systems,
which allows us to determine scaling
limits associated with the infinity limit
of matrix sizes ({\it i.e.} the infinite-particle 
limit) is one of the main topics of
the modern theory of RM \cite{Mehta04}.
The finite and infinite particle systems
showing the orthogonal-unitary and symplectic-unitary
transitions, and transitions between class C
to class CI were studied and 
multitime correlation functions were
determined by Forrester, Nagao and Honner
(FNH) \cite{FNH99}, and by Nagao \cite{N03},
respectively.
The system in the Laguerre ensemble with
$\beta=1$ initial condition reported in the former paper
can be regarded as
the $\nu=\kappa \in \{0, 1,2, \cdots\}$
case of the non-colliding system of the
generalized meanders
and the system reported in the latter paper 
as the $(\nu, \kappa)=(1/2, 1)$ case.

If we think about the system of generalized meanders
apart from the RM theory, however, we can consider the parameters
$\nu$ and $\kappa$ as real numbers, and not necessarily integers nor
half-integers.
In the present paper, we calculate the 
multitime correlation functions of non-colliding
systems of (squared) generalized meanders
for arbitrary values of parameters, 
provided they satisfy the condition
$\nu > -1, \kappa \in [0, 2(\nu+1))$
so that the systems are not collapsed.
We first define the $N$ particle systems 
in a finite time-interval $[0,T]$
and take the $N=T \to \infty$ limit
to construct the two-parameter family
of infinite particle systems.
We prove that the multitime characteristic
functions is given by a Fredholm Pfaffian 
\cite{Rains00} and thus any multitime
correlation function is given by a Pfaffian.
Similarly to the results by FNH \cite{FNH99}
and Nagao \cite{N03} and their 
temporally-homogeneous version (the determinantal
process with the extended Bessel kernel \cite{TW04}),
the elements of the matrix kernels of
Pfaffians are expressed using the Bessel functions,
but we clarify the fact that
they are generally given by the
{\it Riemann-Liouville differintegrals}
of the functions comprising the Bessel functions,
which are used in {\it fractional calculus}
(see, for example, \cite{OS74,Rubin96,Pod99}).
This structure will explain the origin of the
multiple integral expressions for
the elements of the matrix kernels reported by
FNH \cite{FNH99} and Nagao \cite{N03}.

The paper is organized as follows.
In Section 2, the definitions of the generalized
meanders of Yor and their non-colliding systems are
given and the Riemann-Liouville differintegrals 
of the Bessel functions with appropriate factors
are introduced.
The main theorem for the infinite particle limit 
(Theorem \ref{thm:main1}) is then given.
It is demonstrated that,
if we take a further limit in the system
of Theorem \ref{thm:main1}, we will obtain
the temporally homogeneous system of
infinite number of particles,
which is a determinantal process with the extended
Bessel kernel studied in \cite{TW04} 
(see also \cite{Osa03}).
Using the properties of the Riemann-Liouville
differintegrals, we show that
Theorem \ref{thm:main1} includes the results
by FNH \cite{FNH99} and Nagao \cite{N03}
as special cases.
Section 3 is devoted to prove that
for any finite number of particles $N$,
the present system is a
{\it Pfaffian process} (Theorem \ref{thm:finite}), 
in the sense that
any multitime correlation function is given 
by a Pfaffian \cite{Rains00}.
These Pfaffian processes may be regarded as 
the continuous space-time version of the Pfaffian
point processes and Pfaffian Schur processes
studied by Borodin and Rains \cite{BR04}.
Soshnikov used the term
{\it Pfaffian ensembles} in \cite{BS03,Sos03b,Sos04}.
See also \cite{PS02,Joh03,SI04,Ferrari04,IS05}
in the context of study of nonequilibrium phenomena in
the polynuclear growth models,
and \cite{OR03,FS03} in that of shape fluctuations 
of crystal facets.
The processes studied in \cite{FNH99,N03} are
also Pfaffian processes, since the `quaternion
determinantal expressions' of correlation functions, 
introduced and developed by Dyson, Mehta, Forrester, and Nagao 
\cite{D70,Mehta89,Mehta04,NF99,Nag01},
are readily transformed to Pfaffian expressions.
The method of skew-orthogonal functions 
associated with the Laguerre polynomials are
used in Section 4 in order to 
perform matrix inversion and give explicit expressions
for the elements of matrix kernels of Pfaffians.
Asymptotics in $T=N \to \infty$ are studied 
in Section 5.
Appendices are given to show proofs of formulae and
lemmas used in the text.

At the end of this introduction, we would
like to refer to the papers \cite{Ivanov01,CM04},
which reported the further extensions of
RM theory in physics and the representation theory.
We hope that the present paper will demonstrate
the fruitfulness of developing the probability theory of
interacting infinite particle systems
in connection with 
the extensive study of (multi-)matrix models
in the RM theory.

\SSC{Definition of Processes and Results}\label{chap:Results}

\subsection{Non-colliding systems of generalized meanders}
Let $\Z$ and $\R$ be the sets of integers and real numbers,
respectively, and set
$\N = \{1,2,\dots\}$, $\N_0 = \N \cup \{0\}$,
$\Z_- = \Z \setminus \N_0 $, and
$\Rp=\{x \in \R: x \geq 0 \}$.
Let $\Gamma (c), c\in \R\setminus (\Z_- \cup \{0\})$,
be the Gamma function:
$
\Gamma (c )=\int_0^{\infty} dy \ e^{-y}y^{c -1}
$ for $c > 0$,
and
$
\Gamma(c)=
\Gamma(c +[-c]+1)/ \{c(c+1)\cdots (c+[-c]) \}
$
for $c\in (-\infty,0) \setminus \Z_-$,
where $[c]$ is the largest integer 
that is less than or equal to the real number $c$.
For $t > 0$, $x, y \in \Rp$ and $\nu > -1$
we denote by $G_t^{(\nu)}(t;y|x)$
the transition probability density of 
a $2(\nu +1)$-dimensional {\it Bessel process}
\cite{RY98,BS02},
\begin{eqnarray*}
\label{eqn:Bessel1}
G^{(\nu)}(t; y|x)
&=&\frac{y^{\nu+1}}{x^{\nu}}\frac{1}{t} e^{-(x^2+y^2)/2t}
I_{\nu} \left( \frac{x y}{t} \right), 
\quad x>0, y\in\Rp,
\\
G^{(\nu)}(t; y|0)
&=&\frac{y^{2\nu +1}}{2^{\nu} \Gamma (\nu +1)t^{\nu+1}} e^{-y^2 /2t},
\quad y\in\Rp,
\end{eqnarray*}
where $I_{\nu}(z)$ is the modified Bessel function :
$
I_{\nu}(z)= \sum_{n=0}^{\infty} 
(z/2)^{2n+\nu}/\{ \Gamma(n+1)\Gamma(\nu+n+1)\}.
$
For $T>0$, $\kappa \in [0, 2(\nu+1) )$, we put 
\begin{equation*}
h^{(\nu,\kappa)}_{T}(t,x)=\int_0^\infty dy \ G^{(\nu)}(T-t;y|x) 
y^{-\kappa},
\quad x \in\Rp, \, t\in [0,T],
\end{equation*}
and 
\begin{eqnarray}
\label{eqn:GT1}
&&G^{(\nu,\kappa)}_{T}(s,x;t,y)
= \frac{1}{h^{(\nu,\kappa)}_{T}(s,x)}
G^{(\nu)}(t-s;y|x) h^{(\nu,\kappa)}_{T}(t,y),
\quad x > 0, y \in \Rp,
\\
\label{eqn:GT2}
&&G^{(\nu,\kappa)}_{T}(0,0;t,y)
= \frac{\Gamma (\nu+1)}{\Gamma (\nu+1-\kappa/2)}(2T)^{\kappa/2}
G^{(\nu)}(t;y|0) h^{(\nu,\kappa)}_{T}(t,y),
\quad y \in \Rp,
\end{eqnarray}
for $0 \le s \leq t \le T$.
This transition probability density $G_T^{(\nu,\kappa)}(s,x;t,y)$
defines the temporally inhomogeneous 
process in a finite time-interval $[0,T]$,
which is called a {\it generalized meander}.
In particular, when $\nu=1/2$ and $\kappa =1$,
it is identified with the process called a {\it Brownian meander} 
(see Chapter 3 in Yor \cite{Yor92}).

Now we consider the $N$-particle system of generalized meanders
conditioned that they never collide in a finite time-interval
$[0, T]$. 
Let 
$$
\RW = \Big\{\x=(x_{1}, x_{2}, \cdots, x_{N}) 
\in \Rp^{N}:0\leq x_{1} < x_{2} < \cdots < x_{N} \Big\}.
$$
According to the determinantal formula
of Karlin and McGregor \cite{KM59_2},
the transition probability density is given as
\begin{equation}
\g (s, \x; t, \y)
= \frac{f_{N,T}^{(\nu,\kappa)}(s, \x;t,\y)
{\cal N}_{N,T}^{(\nu, \kappa)}(T-t,\y)}
{{\cal N}_{N,T}^{(\nu, \kappa)}(T-s,\x)},
\quad \quad 0\le s \leq t \le T,
\quad \x, \y \in \RW,
\label{eqn:gNnk}
\end{equation}
where
\begin{equation*}
f_{N,T}^{(\nu,\kappa)}(s,\x;t,\y) =
\det_{1 \leq j, k \leq N} \left[G^{(\nu,\kappa)}_{T}(s,x_j,t,y_k)\right],
\quad
{\cal N}_{N,T}^{(\nu, \kappa)}(t,\x)
= \int_{\RW} d \y f_{N,T}^{(\nu,\kappa)}(T-t,\x; T,\y).
\end{equation*}
Since $h_T^{(\nu,0)}(t,x)=1$,
$G_{T}^{(\nu,0)}(s,x;t,y)=G^{(\nu)}(t-s;y|x)$ and thus
$f_{N,T}^{(\nu,0)}$ is temporally homogeneous
and independent of $T$, we will write
$f_{N}^{(\nu)}(t-s ;\y|\x)$ for $f_{N,T}^{(\nu,0)}(s, \x; t,\y)$.
Moreover, note that 
\begin{equation*}
f_{N,T}^{(\nu,\kappa)}(s,\x; t,\y) =
\frac{1}{h^{(\nu,\kappa)}_{T}(s,\x)}f_{N}^{(\nu)}(t-s;\y|\x) 
h^{(\nu,\kappa)}_{T}(t,\y),
\label{eqn:trans}
\end{equation*}
where
$h^{(\nu,\kappa)}_{T}(t,\x) \equiv
\prod_{j=1}^N h^{(\nu,\kappa)}_{T}(t,x_j)$
and
$h^{(\nu, \kappa)}_{T}(T,\x)= \prod_{j=1}^{N} x_{j}^{-\kappa}$.
Then (\ref{eqn:gNnk}) can be written as
\begin{equation}
\g(s,\x;t,\y)
=\frac{1}{\widetilde{{\cal N}}_{N}^{(\nu, \kappa)}(T-s,\x)}
f_{N}^{(\nu)}(t-s;\y|\x)\widetilde{{\cal N}}_{N}^{(\nu, \kappa)}(T-t,\y),
\label{eqn:gNnk1}
\end{equation}
where
\begin{equation}
\widetilde{{\cal N}}_{N}^{(\nu, \kappa)}(t,\x)
= \int_{\RW} d \y \ f_{N}^{(\nu)}(t;\y|\x) \prod_{j=1}^{N} y_{j}^{-\kappa}.
\label{eqn:tildeN}
\end{equation}
In our previous paper \cite{KT04} 
it was shown that, taking the limit 
$\x \to \0 \equiv (0,0,\dots,0)$ at the initial time $s=0$,
(\ref{eqn:gNnk1}) becomes 
\begin{eqnarray}
&&\g(0, {\bf 0}; t, {\bf y}) 
=\CT(t)
\prod_{j=1}^{N} G^{(\nu)}(t, y_{j}|0)
\prod_{1 \leq j < k \leq N} (y_{k}^2-y_{j}^2) \,
\widetilde{{\cal N}}_{N}^{(\nu, \kappa)}(T-t, {\bf y})
\quad
\nonumber\\
\label{eqn:gNnk2}
\end{eqnarray}
for $\nu > -1$ and $\kappa \in [0, 2(\nu+1))$, where
$$
\CT(t)
=\frac{T^{(N+\kappa-1)N/2} t^{-(N-1)N}}{2^{N(N-\kappa-1)/2}}
\prod_{j=1}^{N} \frac{\Gamma(\nu+1) \Gamma(1/2)}
{\Gamma \left(j/2 \right) \Gamma \left( (j+1+2\nu-\kappa)/2 \right)}.
$$
The $N$-particle system of 
{\it non-colliding generalized meanders 
all starting from the origin $\0$ at time $0$}
is defined by the transition probability density $\g$
given above and
it will be denoted by $\X(t) \in \RW, t \in [0,T]$ in the 
present paper.
It makes a two-parameter family of
temporally inhomogeneous processes parameterized by
$\nu > -1$ and $\kappa \in [0, 2(\nu+1))$.

We denote by $\mX$ the space of countable subsets $\xi$ of $\R$
satisfying $\sharp (\xi \cap K) < \infty$
for any compact subset $K$.
For $\x = (x_1, x_2, \dots, x_n)\in \bigcup_{\ell=1}^{\infty} \R^\ell$,
we denote $\{x_i\}_{i=1}^n \in \mX$ simply by $\{\x\}$.
Then $\Xi^{\X}_N (t) = \{\X(t)\}, t \in [0,T]$, is the diffusion process on 
the set $\mX$ with transition density function
$\mg(s, \xi ; t, \eta)$, $0 \leq s \leq t \leq T$:
\begin{equation}
\mg(s, \xi ; t, \eta)=
\left\{
\begin{array}{ll}
\g(s, \x; t, \y),
& \mbox{if} \ s>0,  \ \sharp \xi = \sharp \eta =N,
\\
\g(0, {\bf 0}; t, \y),
& \mbox{if} \ s=0, \ \xi =\{ 0 \}, \ \sharp \eta=N,
\\
0, 
& \mbox{otherwise},
\end{array}\right.
\nonumber
\end{equation}
where $\x$ and $\y$ are the elements of $\RW$ with
$\xi=\{\x\}$, $\eta =\{\y\}$.

For the given time interval $[0,T]$, we consider the $M$
intermediate times
$0 < t_{1} < t_{2} < \cdots < t_{M} < T$.
For convenience, we set
$t_{0}=0$, $t_{M+1}=T$.
For $\x^{(m)} \in \R^N$, $1\leq m \leq M+1$,
and $N'=1,2,\dots, N$, we put 
$\x^{(m)}_{N'} = \left(x_1^{(m)}, x_2^{(m)}, \dots, x_{N'}^{(m)}\right)$
and $\xi_m^{N'} = \{\x^{(m)}_{N'}\}$.
Then the multitime transition density function 
of the process $\Xi^{\X}_N (t)$ is given by
\begin{equation}
\label{def:rho}
\mg \Big(0, \{ 0 \};t_{1}, \xi^N_1;\dots; 
t_{M+1}, \xi^N_{M+1} \Big)
=\prod_{m=0}^{M} \mg \Big(t_{m}, \xi^N_{m}; 
t_{m+1}, \xi^N_{m+1} \Big),
\end{equation}
where we assume $\xi^N_0=\{ 0\}$.
For a sequence $\{N_m \}_{m=1}^{M+1}$ of positive integers 
less than or equal to $N$,
we define the $(N_1, N_2,\dots, N_{M+1})$-multitime 
correlation function by
\begin{eqnarray}
&&\cP_{N,T} \left(t_{1}, \xi_1^{N_{1}}; t_2, \xi_2^{N_{2}};
\dots; t_{M+1}, \xi_{M+1}^{N_{M+1}} \right) 
\nonumber\\
&&=
\int \limits_{\prod_{m=1}^{M+1} \R_+^{N-N_{m}}}
\prod_{m=1}^{M+1}
\frac{1}{(N-N_{m})!}\prod_{j=N_{m}+1}^{N} dx_{j}^{(m)}
\mg \Big(0, \{0\}; t_{1}, \xi^{N}_1; 
\dots;t_{M+1}, \xi^{N}_{M+1} \Big).
\label{def:corr}
\end{eqnarray}

Associated with the generalized meander (\ref{eqn:GT1}),
(\ref{eqn:GT2}), 
we consider a temporally inhomogeneous diffusion process
with transition probability density
\begin{eqnarray*}
\pT(0,0;t,y) &\equiv& G_T^{(\nu, \kappa)}(0,0;t,\sqrt{y}) 
\times \frac{1}{2} y^{-1/2},
\quad y\in\R_{+},
\\
\pT(s,x;t,y) &\equiv& G_T^{(\nu, \kappa)}(s,\sqrt{x}:t, \sqrt{y}) 
\times \frac{1}{2} y^{-1/2},
\quad x>0, y\in\R_{+},
\end{eqnarray*}
$t\in [0,T]$, and call it a {\it squared generalized meander}.
The $N$-particle system of 
{\it non-colliding squared generalized meanders}
$\Y(t), t \in [0,T]$, 
is then defined by
\begin{equation*}
\Y(t) = \Big( X_1(t)^2, X_2(t)^2, \dots, X_N(t)^2 \Big),
\quad t \in [0,T].
\end{equation*}
The correlation function $\cR_{N,T}$ of $\Xi^{\Y}_N(t)=\{\Y(t)\}$
is obtained from (\ref{def:corr}) through the relation
\begin{eqnarray}
&&\cR_{N,T} \left(t_{1}, \zeta_{1}^{N_{1}}; t_2, \zeta_{2}^{N_{2}}; 
\dots; t_{M+1}, \zeta_{M+1}^{N_{M+1}} \right) 
\nonumber\\
&&\qquad =
\cP_{N,T} \left(t_{1}, \xi_{1}^{N_{1}}; t_2, \xi_{2}^{N_{2}}; 
\dots; t_{M+1}, \xi_{M+1}^{N_{M+1}} \right) 
\prod_{m=1}^{M+1}\prod_{j=1}^{N_m}\frac{1}{2x^{(m)}_j},
\label{def:corrY}
\end{eqnarray}
where $\xi_{m}^{N_{m}}=\{\x_{N_{m}}^{(m)} \},
\zeta_{m}^{N_{m}}=\{\y_{N_{m}}^{(m)} \}$ with
$x^{(m)}_j = \sqrt{y^{(m)}_j}$, 
$1\le j\le N_m$, $1\le m \le M+1$.

\subsection{Riemann-Liouville differintegrals of Bessel functions}%

We consider the following 
left and right Riemann-Liouville differintegrals
for integrable functions 
$f$ on $\Rp$,
\begin{eqnarray}
\label{def:RL1}
\D0_{x}^{c} f(x) &=& \frac{1}{\Gamma(n-c)} 
\left( \frac{d}{dx}\right)^n \int_{0}^{x} (x-y)^{n-c-1}
f(y) dy, \\
\label{def:RL2}
_{x}\Dinf^{c} f(x) &=& \frac{1}{\Gamma(n-c)}
\left(- \frac{d}{dx}\right)^n \int_{x}^{\infty}
(y-x)^{n-c-1} f(y) dy,
\end{eqnarray}
where $c \in \R$ and $n=[c+1]_{+}$
with the notation $x_+= \max\{x,0\}$.
It is easy to confirm that, if $c \in \N_{0}$,
both of them are reduced to the ordinary 
multiple derivative,
$$
\D0_{x}^{c} f(x) =(-1)^c  \; _{x}\Dinf^{c} f(x)
= \left( \frac{d}{dx}\right)^c f(x),
$$
and, if $c \in \Z_-$, they are equal to the multiple
integrals,
\begin{eqnarray}
\D0_{x}^{c} f(x) &=& \int_{0}^{x} dy_{|c|-1} 
\int_{0}^{y_{|c|-1}} dy_{|c|-2}
\cdots \int_{0}^{y_{2}} dy_{1} \int_{0}^{y_{1}} d y_{0} f(y_{0}), 
\nonumber\\
_{x}\Dinf^{c} f(x) &=& 
\int_{x}^{\infty} dy_{|c|-1} \int_{y_{|c|-1}}^{\infty} dy_{|c|-2} 
\cdots \int_{y_{2}}^{\infty} dy_{1} \int_{y_{1}}^{\infty} 
dy_{0} f(y_{0}).
\nonumber
\end{eqnarray}
For $c \in (-\infty, 0)\setminus \Z_{-}$
(\ref{def:RL1}) and (\ref{def:RL2}) define
{\it fractional integrals}, and 
for $c \in \Rp \setminus \N_{0}$
{\it fractional differentials}.
The Riemann-Liouville differintegrals are 
most often used in the fractional calculus
(see, for example, \cite{OS74,Rubin96,Pod99}).

Let $J_{\nu}(z)$ be the Bessel functions:
$
J_{\nu}(z)= \sum_{\ell=0}^{\infty}
(-1)^\ell 
(z/2)^{2\ell+\nu}/\{\Gamma(\nu+\ell+1) \ell! \}.
$
We define functions $\Knu$ and $\hKnu$ as
\begin{eqnarray}
&&\Knu (\theta,\eta,x,s)
=(\theta\eta x)^{\nu/2}J_{\nu}(2\sqrt{\theta\eta x})e^{2s\theta\eta}
=e^{2s\theta\eta}
\sum_{\ell=0}^{\infty}
\frac{(-1)^\ell (\theta\eta x)^{\ell+\nu}}{\Gamma(\nu+\ell+1)\ell !},
\label{def:k_nu}
\\
&&\hKnu (\theta,\eta,x,s)
=(\theta\eta x)^{-\nu/2}J_{\nu}(2\sqrt{\theta\eta x})e^{2s\theta\eta}
=e^{2s\theta\eta}
\sum_{\ell=0}^{\infty}
\frac{(-\theta\eta x)^{\ell}}{\Gamma(\nu+\ell+1)\ell !}.
\label{def:hk_nu}
\end{eqnarray}
We will use the following abbreviations
for the Riemann-Liouville differintegrals
of order $c \in \R$
of $\Knu$ and $\hKnu$,
\begin{eqnarray}
\label{eqn:KnuRL}
\Knu^{(c)}(\theta,\eta,x,s)
&=& \D0_{\eta}^{c} \Knu(\theta,\eta,x,s), 
\quad \theta, \eta >0, s \in \R, \\
\label{eqn:hKnuRL}
\hKnu^{(c)}(\theta,\eta,x,s)
&=& _{\eta}\Dinf^{c} \hKnu(\theta,\eta,x,s),
\quad \theta, \eta >0, s < 0.
\end{eqnarray}
We note that, 
if $c \in \R \setminus \N_0$, 
$\Knu^{(c)}$ can be expanded as
\begin{equation}
\Knu^{(c)}(\theta,\eta,x,s)
=\frac{1}{\Gamma (-c)}
\sum_{n=0}^{\infty}\frac{(-1)^n \eta^{n-c}}{n!(n-c)}
\Knu^{(n)}(\theta,\eta,x,s),
\quad \theta,\eta>0, s\in\R.
\label{def:K^c}
\end{equation}
It is also noted that,
since $\hKnu(\theta, \eta,x, s) \to 0$
exponentially fast as $\eta \to \infty$,
if $s \theta < 0$,
\begin{equation}
\hKnu^{(c)}(\theta,\eta,x,s)
= \frac{1}{\Gamma (n-c)}
\int_{\eta}^\infty d \xi \ (\xi-\eta)^{n-c-1}
\hKnu^{(n)}(\theta, \xi, x, s),
\quad \theta, \eta >0, s < 0,
\label{def:hK^c}
\end{equation}
where $n=[c+1]_+$.

\subsection{Results}

We put 
\begin{equation}
\ma =\ma (\nu, \kappa)=\nu - \frac{\kappa}{2},
\qquad
\mb = \mb (\nu, \kappa)=\nu - \kappa,
\label{eqn:aA}
\end{equation}
and introduce functions 
$\cD(s,x;t,y)$, $\cI(s,x;t,y) $, $\cS(s,x;t,y)$ and 
$\tcS (s,x;t,y)$, $x,y \in \Rp$, $s,t < 0$,
\begin{eqnarray}
\cD(s,x;t,y)
&=& - \frac{1}{4(xy)^{\kappa/2}}
\int_0^1 d\theta \ \theta^{1-\kappa}
\Big[\Knu^{(-\mb-1)}(\theta,1,x,-s)\Knu^{(-\mb)}(\theta,1,y,-t)
\nonumber\\
&& \hskip 5cm
-\Knu^{(-\mb)}(\theta,1,x,-s)\Knu^{(-\mb-1)}(\theta,1,y,-t)\Big],
\nonumber\\
\cI(s,x;t,y) 
&=& -(xy)^{\kappa/2}
\int_1^\infty d\theta \ \theta^{\kappa-1}
\Big[\int_1^\infty d\xi \ \xi^\ma 
\hKnu^{(\mb+1)}(\theta,\xi,x,s)\hKnu^{(\mb+1)}(\theta,1,y,t)
\nonumber\\
&& \hskip 5cm
-\hKnu^{(\mb+1)}(\theta,1,x,s)
\int_1^\infty d\xi \ \xi^\ma \hKnu^{(\mb+1)}(\theta,\xi,y,t)\Big],
\nonumber\\
\cS(s,x;t,y)
&=& \frac{1}{2}\left(\frac{x}{y}\right)^{\kappa/2}\int_0^1 d\theta
\Big[ \hKnu^{(\mb+1)}(\theta,1,x,s)\Knu^{(-\mb-1)} (\theta,1,y,-t)
\nonumber\\
&&- 
\big\{\ma\Knu^{(-\mb-1)}(\theta,1,y,-t)
-\Knu^{(-\mb)}(\theta,1,y,-t)\big\}
\int_1^\infty d\xi \ \xi^\ma \hKnu^{(\mb+1)}(\theta,\xi,x,s)
\Big],
\label{eqn:mainDIS}
\end{eqnarray}
and
\begin{eqnarray}
&&\tcS (s,x;t,y)=\cS (s,x;t,y)
- {\bf 1}_{(s<t)} \left(\frac{y}{x}\right)^{\mb/2}\cG(s, x ; t, y),
\label{eqn:mainStilde}
\end{eqnarray}
where
${\bf 1}_{(\omega)}$ is the indicator function:
${\bf 1}_{(\omega)}=1$ if $\omega$ is satisfied
and ${\bf 1}_{(\omega)}=0$ otherwise, and
\begin{equation}
\cG(s, x ; t, y)= 
\int_0^{\infty} d\theta \
J_{\nu}(2\sqrt{\theta x})J_{\nu}(2\sqrt{\theta y})
e^{2(s-t)\theta}.
\label{eqn:mainG}
\end{equation}
For an integer $N$ and a skew-symmetric $2N \times 2N$ matrix
$A=(a_{ij})$, the Pfaffian is defined as
\begin{equation}
\Pf(A) = \Pf_{1 \leq i < j \leq 2N}(a_{ij}) 
= \frac{1}{N !} \sum_{\sigma}
{\rm sgn}(\sigma) a_{\sigma(1) \sigma(2)} a_{\sigma(3) \sigma(4)} 
\cdots a_{\sigma(2N-1) \sigma(2N)},
\label{def:pfaffian}
\end{equation}
where the summation is extended over all permutations $\sigma$
of $(1,2,\dots, 2N)$ with restriction
$\sigma(2k-1) < \sigma(2k), k=1,2,\dots, N$.
We put
$$
\widehat{\Xi}^{\Y}_N(s)
= \{ Y_1(T_N+s), Y_2(T_N+s), \dots, Y_N(T_N+s) \},
\quad s\in [-T_N, 0),
$$
and $\widehat{\Xi}^{\Y}_N(s)=\{ 0\}$, $s\in (-\infty, -T_N)$.
Then we can state the main theorem in the
present paper.
\vskip 3mm
\begin{thm}
\label{thm:main1}
Let $T_N =N$.
Then the process $\widehat{\Xi}^{\Y}_N(s), s\in (-\infty,0)$
converges to the process 
$\widehat{\Xi}^{\Y}_{\infty}(s), s\in (-\infty,0)$,
as $N\to \infty$, 
in the sense of finite dimensional distributions,
whose correlation functions $\cR$ are given by
\begin{equation}
\cR
\left(s_1, \{\y^{(1)}_{N_1}\} ; s_2, \{\y^{(2)}_{N_2}\} ; \dots ; 
s_M, \{\y^{(M)}_{N_M}\} \right)
= \Pf \left[{\cal A}\left(\y^{(1)}_{N_1},\y^{(2)}_{N_2},
\dots, \y^{(M)}_{N_{M}} \right)\right],
\nonumber
\end{equation}
for any $M \geq 1$, 
any sequence $\{N_{m}\}_{m=1}^{M}$ of positive integers,
and any strictly increasing sequence $\{ s_m \}_{m=1}^{M+1}$ of
nonpositive numbers with $s_{M+1}=0$,
where 
${\cal A}\left(\y^{(1)}_{N_1},\y^{(2)}_{N_2},
\dots, \y^{(M)}_{N_{M}} \right)$
is the $2\sum_{m=1}^{M}N_m \times 2\sum_{m=1}^{M}N_m$ 
skew-symmetric matrix defined by 
$$
{\cal A}\left(\y^{(1)}_{N_1},\y^{(2)}_{N_2},
\dots, \y^{(M)}_{N_{M}} \right)=
\Big({\cal A}^{m,n}(y^{(m)}_i, y^{(n)}_j)
\Big)_{1\leq i \leq N_m, 1\leq j \leq N_n,1\leq m,n  \leq M}
$$
with $2\times 2$ matrices ${\cal A}^{m,n}(x,y)$ ;
\begin{eqnarray}
&&{\cal A}^{m,n}(x,y)=
\left( 
\begin{array}{cc}
\cD (s_m,x;s_n,y) & \tcS (s_n,y;s_m,x) 
\cr
-\tcS (s_m,x;s_n,y)& -\cI(s_m,x;s_n,y)
\end{array}
\right).
\nonumber
\end{eqnarray}
\end{thm}

\vskip 3mm

In the infinite-particle system defined by Theorem \ref{thm:main1}, 
we can take the further limit:
$$
s_m \to -\infty \quad \mbox{with the time differences
$s_n -s_m$ fixed}, \quad 1 \leq m,n \leq M.
$$
In this limit, 
$\cD (s_m,x;s_n,y) \cI (s_m,x;s_n,y)\to 0$,
$1 \leq m, n \leq M$,
as we show in Appendix \ref{chap:AppLimit}.
Therefore, we can replace $\cD$ and $\cI$ by
zeros in the matrices. Then the Pfaffian
is reduced to an ordinary determinant of the
$\sum_{m=1}^{M+1} N_{m} \times \sum_{m=1}^{M+1} N_{m}$
matrix,
${\bA}\left(\y^{(1)}_{N_1},\y^{(2)}_{N_2},\dots, 
\y^{(M)}_{N_{M}} \right)
=\left( \ba^{m,n}(y_{i}^{(m)}, y_{j}^{(n)}) 
\right)_{1\leq i \leq N_m, 1\leq j \leq N_n,1\leq m,n  \leq M}$ 
with the elements
$$
\ba^{m,n}\left(y^{(m)}_i,y^{(n)}_j\right)
=\widetilde{\bS}\left(s_m,y^{(m)}_i;s_n,y^{(n)}_j\right),
$$
where
\begin{eqnarray}
&&\widetilde{\bS}(s, x ; t, y)=
\left\{
   \begin{array}{ll}
\displaystyle{
\int_0^1 d\theta \
J_{\nu}(2\sqrt{\theta x})J_{\nu}(2\sqrt{\theta y})
e^{2(s-t)\theta},}
 & \mbox{if} \ s > t, \\
 & \\
\displaystyle{
\frac{J_{\nu}(2\sqrt{x})\sqrt{y}J_{\nu}'(2\sqrt{y})
-J_{\nu}(2\sqrt{y})\sqrt{x}J_{\nu}'(2\sqrt{x})}{x-y }},
 & \mbox{if} \ s=t, \\
 & \\
\displaystyle{
-\int_1^{\infty} d\theta \
J_{\nu}(2\sqrt{\theta x})J_{\nu}(2\sqrt{\theta y})
e^{2(s-t)\theta},}
 & \mbox{if} \ s< t, \\
   \end{array}\right. 
\nonumber
\end{eqnarray}
with $J_{\nu}'=d J_{\nu}(z)/dz$.
Hence, in this limit we obtain a temporally 
homogeneous system of infinite number of particles, 
whose correlation functions are given by
\begin{equation}
\widetilde{\rho}^{\Y}
\left(s_1, \{\y^{(1)}_{N_1} \} ; s_2, 
\{\y^{(2)}_{N_2} \} ; \dots ; 
s_M, \{\y^{(M)}_{N_M}\} \right)
= \det {\bA} \left(\y^{(1)}_{N_1},\y^{(2)}_{N_2},
\dots, \y^{(M)}_{N_{M}} \right).
\label{eqn:rhodash}
\end{equation}

\vskip 0.3cm
\noindent{\bf Remark 1.} 
Forrester, Nagao, and Honner \cite{FNH99} studied
the orthogonal-unitary and symplectic-unitary universality
transitions in random matrix theory by
giving the quaternion determinantal expressions
of (two-time) correlation functions for 
parametric RM models.
One of their results for the `Laguerre ensemble
with $\beta=1$ initial condition', which shows
the orthogonal-unitary transition, can be
reproduced from Theorem \ref{thm:main1}
by setting
$$
{\rm (i)} \qquad \kappa=\nu \qquad
\Longleftrightarrow \qquad
\ma = \frac{\nu}{2}, \quad \mb=0, \qquad
\mbox{where} \quad \nu \in \N_0.
$$
This fact may be readily seen, if we notice that
by definition
\begin{eqnarray}
\widetilde{J}_{\nu}^{(-1)} (\theta,1,x, s) 
&=& \int_{0}^{1} d \eta \, (\theta \eta x)^{\nu/2}
J_{\nu}(2 \sqrt{\theta \eta x}) e^{2 s \theta \eta}
\nonumber\\
&=& \theta^{-1} x^{\nu/2}
\int_{0}^{\theta} du \, u^{\nu/2} 
J_{\nu}(2 \sqrt{ux}) e^{2 s u}.
\nonumber
\end{eqnarray}

\vskip 0.3cm
\noindent{\bf Remark 2.}
Nagao's result on the multitime correlation functions
for {\it vicious random walk with a wall} \cite{N03} can be
regarded as the special case of Theorem \ref{thm:main1},
in which
$$
{\rm (ii)} \qquad \nu=\frac{1}{2}, \quad \kappa=1
\qquad 
\Longleftrightarrow \qquad
\ma = 0, \quad \mb=-\frac{1}{2}.
$$
This fact can be confirmed by noting that,
by definition (\ref{eqn:KnuRL}) with 
$\widetilde{J}_{1/2}(\theta, 0; x, s)=0$,
\begin{eqnarray}
\widetilde{J}^{(-1/2)}_{1/2}(\theta, 1, x, s)
&=& \frac{(\theta x)^{1/4}}{\sqrt{\pi}} \int_{0}^{1} d \eta \,
(1-\eta)^{-1/2} \eta^{1/4} 
J_{1/2}(2 \sqrt{ \theta \eta x}) e^{2 s \theta \eta}, \nonumber\\
\widetilde{J}^{(1/2)}_{1/2}(\theta, 1, x, s)
&=& \frac{(\theta x)^{1/4}}{\sqrt{\pi}} \int_{0}^{1} d \eta \,
(1-\eta)^{-1/2} \frac{d}{d \eta} \left\{
\eta^{1/4} J_{1/2}(2 \sqrt{ \theta \eta x}) e^{2 s \theta \eta}
\right\}, \nonumber
\end{eqnarray}
by (\ref{def:hK^c}),
$$
\widehat{J}_{1/2}^{(1/2)}(\theta, \eta, x, s )
= - \frac{(\theta x)^{-1/4}}{\sqrt{\pi}}
\int_{\eta}^{\infty} d \xi \,
(\xi-\eta)^{-1/2}
\frac{d}{d\xi} \left\{
\xi^{-1/4} J_{1/2} (2 \sqrt{\theta \xi x})
e^{2 s \theta \xi} \right\}, \quad s < 0,
$$
and, by definition (\ref{eqn:hKnuRL}),
$$
\int_{1}^{\infty} d \xi \widehat{J}_{1/2}^{(1/2)}(\theta, \xi, x, s)
= \frac{(\theta x)^{-1/4}}{\sqrt{\pi}}
\int_{1}^{\infty} d \eta \,
(\eta-1)^{-1/2} \eta^{-1/4}
J_{1/2}(2 \sqrt{\theta \eta x}) e^{2 s \theta \eta},
\quad s < 0.
$$
In this case, the system shows the transition
between the class C and class CI of 
the Bogoliubov-de Gennes universality classes of
nonstandard RM theory \cite{N03,KTNK03,KT04}.

\vskip 0.3cm
\noindent{\bf Remark 3.}
From the results for finite 
non-colliding processes \cite{KT04},
we expect that, when
$$
{\rm (iii)} \qquad \kappa=\nu+1 \qquad
\Longleftrightarrow \qquad
\ma = \frac{\nu-1}{2}, \quad \mb=-1, \qquad
\mbox{where} \quad \nu \in \N_0,
$$
the present infinite particle system will show the
transition from the chiral GUE to the 
chiral GOE  of the universality classes
and when 
$$
{\rm (iv)} \qquad \nu=-\frac{1}{2}, \quad \kappa=0
\qquad 
\Longleftrightarrow \qquad
\ma = \mb=-\frac{1}{2},
$$
that from the class D to the `real-component version'
of class D of the 
Bogoliubov-de Gennes universality classes \cite{KT04}.

\vskip 0.3cm
\noindent{\bf Remark 4.}
Following the argument given in
\cite{PS02,Joh03}, tightness in time
can be proved and transition phenomena
observed in the limit $s_{M} \to 0$ 
may be generally discussed, 
which will be reported elsewhere.

\vskip 0.3cm
\noindent{\bf Remark 5.}
The homogeneous system (\ref{eqn:rhodash}) 
was studied in \cite{TW04,Osa03}.


\SSC{Correlation Functions Given by Pfaffians}\label{chap:Correlation}

\subsection{The multitime transition density}

If we put 
\begin{eqnarray*}
&&\widetilde{G}^{(\nu, \kappa)}(t, y|x)
= G^{(\nu)}(t, y|x) \times
\left( \frac{y}{x} \right)^{-\kappa},
\quad x>0, y \in \R_+,
\label{eqn:tildep}
\\
&&\widetilde{G}^{(\nu, \kappa)}(t, y|0)
= G^{(\nu)}(t, y|0) \times y^{-\kappa},
\quad y \in \R_+,
\label{eqn:tildep2}
\end{eqnarray*}
the multitime transition density 
(\ref{def:rho}) with $t_{0}=0, t_{M+1}=T$,
and $\xi_{0}=\{0\}$ is written as
\begin{eqnarray*}
&&\mg \Big(0,\{ 0 \};t_{1},\{\x_N^{(1)}\};
\cdots; t_{M+1},\{\x_N^{(M+1)}\} \Big) 
\nonumber\\
&&\qquad = \CT(t_1)
\prod_{1 \leq j < k \leq N} \left\{(x_{k}^{(1)})^2-(x_{j}^{(1)})^2 \right\}
\prod_{1 \leq j < k \leq N} {\rm sgn}(x_{k}^{(M+1)}-x_{j}^{(M+1)})
\nonumber\\
&& \qquad\qquad \times 
\prod_{j=1}^{N} \widetilde{G}^{(\nu, \kappa)}(t_1,x_{j}^{(1)} |0)
\prod_{m=1}^{M}\det_{1 \leq j, k \leq N} \Bigg[
\widetilde{G}^{(\nu, \kappa)}
(t_{m+1}-t_{m}, x_{j}^{(m+1)}|x_{k}^{(m)}) \Bigg],
\label{eqn:g2b}
\end{eqnarray*}
where (\ref{eqn:gNnk1}) and (\ref{eqn:gNnk2})
with (\ref{eqn:tildeN}) are used.

Through the relation (\ref{def:corrY}), 
the multitime transition density for the process
$\{ \Y(t) \}, t \in [0,T]$, denoted by $\mp$ is then written as
\begin{eqnarray}
&&\mp \Big(0, \{ 0\}; t_{1}, \{\y^{(1)}\}; \dots;
t_{M+1}, \{\y^{(M+1)}\} \Big)
\nonumber\\
&&\qquad=\CT(t_1) h_{N}(\y^{(1)}){\rm sgn} 
\left( h_{N}(\y^{(M+1)}) \right)
\nonumber\\
&&\qquad\qquad\times \prod_{k=1}^N \p(t_1,y_k^{(1)}|0)
\prod_{m=1}^{M} \det_{1 \leq j, k \leq N}
\Bigg[ \p(t_{m+1}-t_m, y_{j}^{(m+1)} | y_{k}^{(m)}) \Bigg],
\label{eqn:gNT2}
\end{eqnarray}
where
$$
h_N(\y) \equiv \prod_{1\le i< j \le N}(y_j - y_i),
\quad \y\in\R^N,
$$
\begin{eqnarray}
\p(t,y |0) &\equiv& \widetilde{G}^{(\nu, \kappa)}(t, \sqrt{y}|0) 
\times \frac{1}{2} y^{-1/2}
\nonumber\\
&=&
\frac{y^{\ma}}{2^{\nu+1}\Gamma(\nu+1)t^{\nu+1}} e^{-y/2t},
\quad y\in\R_+,
\nonumber
\end{eqnarray}
\begin{eqnarray}
\p(t-s,y|x) &\equiv& \widetilde{G}^{(\nu, \kappa)}(t-s, \sqrt{y}|\sqrt{x}) 
\times \frac{1}{2} y^{-1/2} 
\nonumber\\
&=& 
\frac{e^{-(x+y)/\{2(t-s)\}}}{2(t-s)}
\left(\frac{y}{x}\right)^{\mb/2} 
I_{\nu} \left( \frac{\sqrt{xy}}{t-s} \right),
\ x>0, \ y\in\R_+.
\label{eqn:gmxy}
\end{eqnarray}
Expectations related to the process
$\{\Y(t_1)\}, \{\Y(t_2)\},\dots, \{\Y(t_{M+1})\}$ 
are denote by $\E_{N,T}^{\Y}$ :
\begin{eqnarray}
&&\E_{N,T}^{\Y} 
\Big[ f(\{\Y(t_1)\}, \{\Y(t_2)\},\dots, \{\Y(t_{M+1})\}) \Big]
= \left( \frac{1}{N!} \right)^{M+1} \int_{\Rp^{N(M+1)}}
\prod_{m=1}^{M+1} d \y^{(m)} \,
\nonumber\\
&&\qquad\qquad \times
f(\{\y^{(1)} \}, \{\y^{(2)}\},\dots, \{\y^{(M+1)}\}) \,
\mp \Big(0, \{ 0 \}; t_{1}, \{\y^{(1)}\}; 
\dots; t_{M+1}, \{\y^{(M+1)}\} \Big).
\label{eqn:expect}
\end{eqnarray}

\subsection{Fredholm Pfaffian representation
of characteristic function
and Pfaffian process}\label{section32}

For simplicity of expressions, we assume from now on
that the number of particles $N$ is even.
The references \cite{NF99,N03} will be useful to give necessary
modifications to the following expressions in the 
case that $N$ is odd.
Let $C_{0}(\R)$ be the set of all continuous real functions
with compact supports. 
For ${\bf f}=(f_{1}, f_{2}, \cdots, f_{M+1}) \in C_{0}(\R)^{M+1}$,
and $\vtheta=(\theta_{1}, \theta_{2}, \cdots, \theta_{M+1})
\in \R^{M+1}$,
the multitime characteristic function is defined
for the process $\{\Y(t)\}, t \in [0,T]$ as
\begin{equation}
\kPsi ( {\bf f}; \vtheta)
= \E_{N,T}^{\Y} \left[
\exp \left\{ \sqrt{-1} \sum_{m=1}^{M+1}
\theta_{m} \sum_{i_{m}=1}^{N} 
f_{m}(Y_{i_{m}}(t_{m})) \right\} \right]
\label{eqn:Psi0}
\end{equation}
Let 
$
\chi_{m}(x)=e^{\sqrt{-1} \theta_{m} f_{m}(x)}-1, 
\, 1 \leq m \leq M+1.
$
Then by the definition of
multitime correlation function (\ref{def:corrY})
with (\ref{def:corr}), we have 
\begin{eqnarray}
\kPsi ( {\bf f}; \vtheta )
&=& \sum_{N_{1}=0}^{N} \sum_{N_{2}=0}^{N} \cdots
\sum_{N_{M+1}=0}^{N}
\prod_{m=1}^{M+1}\frac{1}{N_m !}
\int_{\Rp^{N_{1}}} d \y_{N_{1}}^{(1)}
\int_{\Rp^{N_{2}}} d \y_{N_{2}}^{(2)} \cdots
\int_{\Rp^{N_{M+1}}} d \y_{N_{M+1}}^{(M+1)} \nonumber\\
&& \times \prod_{m=1}^{M+1} \prod_{i^{(m)}=1}^{N_{m}} 
\chi_{m} \Big(y_{i^{(m)}}^{(m)} \Big)
\cR_{N,T} \left(t_{1}, \{\y_{N_{1}}^{(1)} \}; 
t_{2}, \{\y_{N_{2}}^{(2)} \};
\dots ; t_{M+1}, \{\y_{N_{M+1}}^{(M+1)}\} \right),
\qquad
\label{eqn:kPhiA1}
\end{eqnarray}
that is, the multitime characteristic function
is a generating function of multitime correlation
functions $\rho_{N,T}^{\Y}$.

We consider a vector space ${\cal V}$ with the orthonormal
basis $\Big\{ |m,x \rangle \Big\}_{1 \leq m \leq M+1,
x \in \Rp}$, which satisfies
\begin{equation}
 \langle m, x | n, y \rangle 
 =\delta_{m n} \delta(x-y),
 \quad m, n =1,2, \cdots, M+1, x, y \in \Rp,
\label{eqn:orthZ}
\end{equation}
where $\delta_{m n}$ and $\delta(x-y)$
denote Kronecker's delta and
Dirac's $\delta$-measure, respectively.
We introduce the operators $\hat{J}, \hat{p},
\hat{p}_{+}, \hat{p}_{-}$ and $\hat{\chi}$ 
acting on ${\cal V}$ as follows
\begin{eqnarray}
\label{eqn:hatJ}
\langle m, x | \hat{J} | n, y \rangle 
&=& {\bf 1}_{(m=n=M+1)} {\rm sgn}(y-x) , \\
\langle m, x | \hat{p} | n, y \rangle 
&=& {\bf 1}_{(m < n)} \p(t_{n}-t_{m}, y|x) 
+ {\bf 1}_{(m > n)} \p(t_{m}-t_{n}, x | y) \nonumber\\
\label{eqn:hatp}
&+& {\bf 1}_{(m=n)} \delta(x-y), \\
\label{eqn:p+-}
\langle m, x | \hat{p}_{+} | n, y \rangle 
&=& {\bf 1}_{(m < n)} \p(t_{n}-t_{m}, y | x)
= \langle n, y| \hat{p}_{-} | m, x \rangle, \\
\label{eqn:hatchi}
\langle m, x | \hat{\chi} | n, y \rangle
&=& \chi_{m}(x) \delta_{m n} \delta(x-y),
\end{eqnarray}
and we will use the convention
$$
\langle m, x | \hat{A} | n, y \rangle
\langle n, y | \hat{B} | \ell, z)
=\sum_{n=1}^{M+1} \int_{\Rp} dy \,
A(m,x; n,y)B(n,y;\ell,z)
=\langle m,x | \hat{A}\hat{B}|\ell,z \rangle
$$
for operators $\hat{A}$ and $\hat{B}$ with
$\langle m, x | \hat{A} | n, y \rangle=A(m,x;n,y)$ and
$\langle m, x | \hat{B} | n, y \rangle=B(m,x;n,y)$.

Let $M_{i}(x)$ be an arbitrary polynomial
of $x$ with degree $i$ in the form
$M_{i}(x)=b_{i} x^{i}+ \cdots$ 
with a constant $b_{i} \not=0$
for $i \in \N_{0}$. Since the product of differences
$h_{N}(\x)$ is equal to the Vandermonde determinant, we have
\begin{equation}
h_{N}(\x)=\left\{ \prod_{k=1}^{N} b_{k-1} \right\}^{-1}
\det_{1 \leq i, j \leq N} \Big[
M_{i-1}(x_{j}) \Big].
\label{eqn:hN2}
\end{equation} 
Then we consider the set of linearly independent vectors 
$\Big\{ |i \rangle \, ; \,i \in \N \Big\}$ in ${\cal V}$
defined by
$$
|i\rangle = |m,x \rangle \langle m,x|i\rangle,
$$
where
\begin{equation}
\langle m, x | i \rangle
=\langle i | m, x \rangle 
=\int_{\Rp} dy M_{i-1}(y) \p(t_{1}, y | 0 )
\p(t_{m}-t_{1}, x | y),
\label{eqn:M}
\end{equation}
$i \in \N, m=1,2, \dots, M+1, x \in \Rp$.
We will use the convention
$$
\langle i | \hat{A} | j \rangle
\langle j | \hat{B} | m, x \rangle
= \sum_{j=1}^{\infty} A_{ij} B_{j}^{(m)}(x)
=\langle i| \hat{A}\circ \hat{B} | m,x \rangle,
$$
for
$A_{ij} = \langle i | \hat{A} | j \rangle$
and
$B_{j}^{(m)}(x) =\langle j | \hat{B} | m, x \rangle.$
It should be noted that the vecors 
$\Big\{ |i \rangle \, ; \, i \in \N \Big\}$
are not assumed
to be mutually orthogonal. 
By these vectors, however, any operator $\hat{A}$ on ${\cal V}$
may have a semi-infinite matrix representation
$A=\Big(\langle i | \hat{A} | j \rangle 
\Big)_{i,j \in \N}$.
If the matrix $A$ representing an operator $\hat{A}$
is invertible, we define the operator
$\hat{A}^{\bigtriangleup}$ so that its matrix representation
is the inverse of $A$;
\begin{equation}
\label{eqn:Delta1}
  \Big( \langle i | \hat{A}^{\bigtriangleup} | j \rangle
  \Big)_{i,j \in \N} = A^{-1},
\end{equation}
that is,
$
\langle i | \hat{A} |j \rangle
\langle j | \hat{A}^{\bigtriangleup} | k \rangle
=\langle i| \hat{A} \circ \hat{A}^{\bigtriangleup} 
| k \rangle = \delta_{ik}, \, i,k \in \N.
$

Let ${\cal P}_{N}$ be a linear operator projecting
${\it Span}\Big\{ |i \rangle \, ; \, i \in \N \Big\}$ 
to its $N$-dimensional subspace
${\it Span}\Big\{ |i \rangle \, ; \, i=1,2, \dots, N \Big\}$ 
such that
$$
\langle i | {\cal P}_{N} | m, x \rangle
= \langle m, x | {\cal P}_{N} | i \rangle
= \left\{
\begin{array}{ll}
\langle i |m, x \rangle,
& \mbox{if} \ 1 \leq i \leq N,
\\
& \\
0, 
& \mbox{otherwise}.
\end{array}\right.
$$
We will use the abbreviation
$\hat{A}_{N}={\cal P}_{N} \hat{A} {\cal P}_{N}$
for an operator $\hat{A}$.
If the $N \times N$ matrix defined by
$A_N=( \langle i | \hat{A}_N | j \rangle)_{1 \leq i, j \leq N}$
is invertible, then $(\hat{A}_N)^{\bigtriangleup}$ is defined
so that 
$\Big(\langle i | (\hat{A}_N)^{\bigtriangleup} | j \rangle \Big)
_{1 \leq i, j \leq N}=(A_N)^{-1}$,
and $\langle i | (\hat{A}_N)^{\bigtriangleup} | j \rangle =0$,
if $i \ge N+1$ or $j \ge N+1$.

As shown in \ref{chap:Proof(3.16)}, we can prove that 
\begin{equation}
\Bigg\{ \kPsi ( {\bf f}; \vtheta) \Bigg\}^2
= {\rm Det} \left( I_{2} \delta_{m n} \delta(x-y)
+ \left( 
\begin{array}{cc}
\tS^{m,n}(x,y) & \tI^{m,n}(x,y) \cr
D^{m,n}(x,y) & \tS^{n,m}(y,x)
\end{array}
\right)
\chi_{n}(y) \right),
\label{eqn:Psi2}
\end{equation}
where ${\rm Det}$ denotes the Fredholm determinant.
Here $I_{2}$ is the unit matrix with size 2,
\begin{eqnarray}
D^{m,n}(x,y) 
&=& -\langle m,x | 
\circ (\hat{J}_{N})^{\bigtriangleup} \circ 
| n, y \rangle, 
\nonumber\\
S^{m,n}(x,y) 
&=& \langle m, x | \hat{p} \hat{J}
\circ (\hat{J}_{N})^{\bigtriangleup} \circ 
| n, y \rangle, 
\nonumber\\
I^{m,n}(x,y)
&=& -\langle m, x | \hat{p} \hat{J} 
\circ (\hat{J}_{N})^{\bigtriangleup} \circ 
\hat{J} \hat{p} | n, y \rangle,
\label{eqn:DSI1}
\end{eqnarray}
and
\begin{eqnarray}
\tS^{m,n}(x,y) &=& S^{m,n}(x,y)
-\langle m,x | \hat{p}_{+} | n, y \rangle \nonumber\\
\tI^{m,n}(x,y) &=& I^{m,n}(x,y)
+ \langle m, x | \hat{p} \hat{J} \hat{p} | n, y 
\rangle. 
\label{eqn:SItilde1}
\end{eqnarray}
It implies that the multitime characteristic function
is given by the {\it Fredholm Pfaffian} \cite{Rains00},
\begin{equation}
\label{eqn:Psi3}
\kPsi( {\bf f}; \vtheta)
={\rm PF} \Big(
J_{2} \delta_{m n} \delta(x-y) 
+ \sqrt{\chi_{m}(x)} A^{m,n}(x,y)
\sqrt{\chi_{n}(y)} \Big),
\end{equation}
where $ J_{2}=\displaystyle{\left( 
\begin{array}{cc}
 0 & 1 \cr -1 & 0
\end{array}
\right)}$ and
\begin{eqnarray}
A^{m,n}(x,y) 
&=& J_{2} \left(
\begin{array}{cc}
  \tS^{m,n}(x,y) & \tI^{m,n}(x,y) \cr
  D^{m,n}(x,y) & \tS^{n,m}(y,x)
\end{array}
 \right) \nonumber\\
&=&\left(
\begin{array}{cc}
 D^{m,n}(x,y) & 
\tS^{n,m}(y,x) \cr
-\tS^{m,n}(x,y) & -\tI^{m,n}(x,y) 
\end{array}
 \right).
\label{eqn:matrixA}
\end{eqnarray}
It is defined by
\begin{eqnarray}
&& {\rm PF} \Big(
J_{2} \delta_{m n} \delta(x-y) 
+ \sqrt{\chi_{m}(x)} A^{m,n}(x,y)
\sqrt{\chi_{n}(y)} \Big) \nonumber\\
&& \qquad = \sum_{N_{1}=0}^{N} \sum_{N_{2}=0}^{N} \cdots
\sum_{N_{M+1}=0}^{N}
\prod_{m=1}^{M+1}\frac{1}{N_m !}
\int_{\Rp^{N_{1}}} d \y_{N_{1}}^{(1)}
\int_{\Rp^{N_{2}}} d \y_{N_{2}}^{(2)} \cdots
\int_{\Rp^{N_{M+1}}} d \y_{N_{M+1}}^{(M+1)} \nonumber\\
\label{eqn:corr3}
&& \qquad \qquad \qquad 
\prod_{m=1}^{M+1} \prod_{i^{(m)}=1}^{N_{m}} 
\chi_{m} \Big(y_{i^{(m)}}^{(m)} \Big) {\rm Pf} \Bigg(
A \Big( \y_{N_{1}}^{(1)}, \y_{N_{2}}^{(2)}, \cdots,
\y_{N_{M+1}}^{(M+1)} \Big) \Bigg),
\label{eqn:kPhiA2}
\end{eqnarray}
where 
$A \Big( \y_{N_{1}}^{(1)}, \y_{N_{2}}^{(2)}, \cdots,
\y_{N_{M+1}}^{(M+1)} \Big)$ denotes the 
$2 \sum_{m=1}^{M+1} N_{m} \times 2 \sum_{m=1}^{M+1} N_{m}$
skew-symmetric matrices constructed
from (\ref{eqn:matrixA}) as
$$
A \Big( \y_{N_{1}}^{(1)}, \y_{N_{2}}^{(2)}, \cdots,
\y_{N_{M+1}}^{(M+1)} \Big)
= \left( A^{m,n} (y_{i}^{(m)}, y_{j}^{(n)} )
\right)_{1 \leq i \leq N_{m}, 1 \leq j \leq N_{n},
1 \leq m, n \leq M+1}
$$
for $N_{m}=1,2, \cdots, N, 1 \leq m \leq M+1$.
Comparison of (\ref{eqn:kPhiA1}) and 
(\ref{eqn:Psi3}) with (\ref{eqn:kPhiA2})
immediately gives the following statement.
\begin{thm}
\label{thm:finite}
The $N$-particle non-colliding system of
squared generalized meanders 
$\Y(t), t \in [0,T]$ is a 
Pfaffian process, in the sense that
any multitime correlation function is given by a Pfaffian
\begin{equation*}
\label{eqn:rho2}
\cR_{N,T} 
\left(t_{1}, \{\y^{(1)}_{N_1}\}; t_2, \{\y^{(2)}_{N_2}\}; 
\dots; t_{M+1}, \{\y^{(M+1)}_{N_{M+1}}\} \right) 
={\rm Pf} \left(
A \Big( \y_{N_{1}}^{(1)}, \y_{N_{2}}^{(2)}, \cdots,
\y_{N_{M+1}}^{(M+1)} \Big) \right).
\end{equation*}
\end{thm}
\vskip 0.3cm

\SSC{Skew-Orthogonal Functions and Matrix Inversion }\label{chap:Skew}

\subsection{Skew-symmetric inner products}

Consider the $N \times N$ skew-symmetric matrix
$A_{0} =( (A_{0})_{ij})_{1 \leq i,j \leq N}$
with
\begin{equation}
(A_{0})_{ij} = \langle i | \hat{J}_{N} | j \rangle
=\langle i | m, x \rangle \langle m, x |
\hat{J} | n, y \rangle \langle n, y | j \rangle,
\quad i,j=1,2,\dots,N.
\label{eqn:A0B}
\end{equation}
In order to clarify the fact that each element
$(A_{0})_{ij}$ is a functional of the polynomials
$M_{i-1}(x)$ and $M_{j-1}(x)$ through (\ref{eqn:M}),
we introduce the skew-symmetric inner product
\begin{equation}
\langle f, g \rangle \equiv
\int_0^\infty  dx \int_0^\infty dy \ F(x,y)
\p(t_1,x | 0)\p(t_1,y | 0) f(x) g(y),
\label{eqn:inner1}
\end{equation}
where
\begin{equation}
F(x,y)
= \int_{0}^{\infty} dw \int_{0}^{w} dz \ 
\left| \begin{array}{cc}
 \p(T-t_1, z|x) & \p(T-t_1, w|x)
 \cr
 \p(T-t_1,z|y) & \p(T-t_1,w|y)
\end{array}
\right|, \quad x, y \in \Rp.
\label{def:F}
\end{equation}
Then we have the expression
\begin{equation}
(A_{0})_{ij}=\langle M_{i-1}, M_{j-1} \rangle, \quad
i,j =1,2, \cdots, N.
\label{eqn:A0B2}
\end{equation}

We now rewrite the skew-symmetric inner product 
(\ref{eqn:inner1})
by using the simpler one 
\begin{eqnarray}
\langle f, g \rangle_{*} &=& - \langle g, f \rangle_{*}
\nonumber\\
&\equiv& \int_{0}^{\infty} dw \, e^{-w/2} w^{\ma}
\int_{0}^{w} dz \, e^{-z/2} z^{\ma}
\Bigg\{ f(z) g(w)-f(w) g(z) \Bigg\},
\label{eqn:eskew1}
\end{eqnarray}
which we call
the {\it elementary skew-symmetric inner product.}
Remind that $\p$ is given by (\ref{eqn:gmxy}) using
the modified Bessel function. We will
expand it in terms of the Laguerre polynomials,
$L_j^\alpha (x)=(x^{-\alpha}e^x/j!)
(d/dx)^j(e^{-x}x^{j+\alpha})$, 
$\alpha \in \R$, $j \in \N_{0}$,
using the formula
\begin{equation}
\sum_{j=0}^\infty
\frac{\Gamma(j+1)L^{\nu}_j(x)L^{\nu}_j(y)r^j}{\Gamma(j+1+\nu)}
= \frac{1}{1-r}e^{-\frac{(x+y)r}{1-r}}(xyr)^{-\nu/2}
I_\nu\left(\frac{2\sqrt{xyr}}{1-r}\right), \quad
|r| < 1, \nu > -1.
\label{eqn:Mehler}
\end{equation}
(See the corresponding calculation for the
non-colliding Brownian particles in \cite{KNT04},
where the heat kernel was expanded in terms
of the Hermite polynomials.)
For this purpose, it is useful to introduce the variables
$$
c_n=\frac{t_n(2T-t_n)}{T}, 
\qquad \chi_n =\frac{2T-t_n}{t_n}, \quad
n=1,2, \cdots, M+1,
$$
since we can see that
\begin{eqnarray*}
&&\p (t_{n}-t_m, c_{n} \eta |c_m \xi)
= \frac{1}{2(t_n-t_m)}
I_{\nu}\left(\frac{2\sqrt{\xi\eta\chi_n/\chi_m}}
{1-\chi_n/\chi_m}\right)
\left( \frac{c_n\eta}{c_m\xi} \right)^{\mb/2}
\\
&&\qquad\quad \times \exp\left[
-\left( \frac{1}{1-\chi_n/\chi_m}-1+\frac{t_m}{2T}\right)\xi
-\left( \frac{1}{1-\chi_n/\chi_m}-\frac{t_n}{2T}\right)\eta
\right],
\end{eqnarray*}
and, if we apply the formula 
(\ref{eqn:Mehler}) with 
$r=\chi_{n}/\chi_{m}, x=\xi$ and $y=\eta$, 
it is written as 
\begin{eqnarray}
\p (t_{n}- t_m, c_{n} \eta | c_m\xi)
&=&\left( \frac{t_m}{t_n} \right)^{\nu+1}
c_m^{-\ma-1}\xi^{\kappa/2} (c_n \eta)^{\ma}
\exp\left[-\frac{t_m}{2T}\xi-(1- \frac{t_n}{2T})\eta\right]
\nonumber\\
&& \times 
\sum_{j=0}^\infty\frac{\Gamma(j+1)}{\Gamma(j+1+\nu)}
\left(\frac{\chi_n}{\chi_m}\right)^j L^{\nu}_j(\xi)L^{\nu}_j(\eta).
\label{eqn:expansionL}
\end{eqnarray}
That is, $c_n$ and $\chi_n$ give 
the spatial scale of spread of $N$ particles
and the proper temporal factor at time $t_n$,
respectively.
(See equation (17) and explanation below it
in \cite{NKT03}, where the variable $c_n$ was
determined by showing that the one-particle
density obeys Wigner's semicircle law scaled by
$c_n$ for the non-colliding Brownian particles.)
In particular, for $n=M+1$ we have
\begin{eqnarray}
\p (T-t_m, T \eta |c_m\xi)
&=&\frac{t_m^{\nu+1}}{T^{\kappa/2+1}}
c_m^{-\ma-1}\xi^{\kappa/2} \eta^{\ma}
\exp\left[\left(1- \frac{t_m}{T}\right)\frac{\xi}{2}\right]
\nonumber\\
&& \times
e^{-\xi/2}e^{-\eta/2}
\sum_{j=0}^\infty\frac{\Gamma(j+1)}{\Gamma(j+1+\nu)}
\chi_m^{-j} L^{\nu}_j(\xi)L^{\nu}_j(\eta),
\label{eqn:gmM+1}
\end{eqnarray}
since $c_{M+1}=T$ and $\chi_{M+1}=1$.
Then we obtain the relation
\begin{eqnarray}
\left\langle f\left(\frac{\cdot}{c_1}\right), 
g\left(\frac{\cdot}{c_1}\right)\right\rangle
&=& \frac{2^{-2\nu-2}T^{-\kappa}}{\Gamma(\nu+1)^2}
\int_0^{\infty}dx\int_0^{\infty}dy \, 
e^{-x}e^{-y}x^{\nu}y^{\nu}f(x)g(y)
\nonumber\\
&\times&
\sum_{j=0}^{\infty}\sum_{k=0}^{\infty}\chi_1^{-j-k}
L^{\nu}_j(x)L^{\nu}_k(y)
\Big\langle \frac{\Gamma(j+1)}{\Gamma(j+1+\nu)}L_j^{\nu},
\frac{\Gamma(k+1)}{\Gamma(k+1+\nu)}L_k^{\nu}\Big\rangle_*.
\label{eqn:fg_1}
\end{eqnarray}

\subsection{Skew-orthogonal polynomials}

For $\alpha \in \R$ and $n\in \Z$ we define
\begin{equation}
{n+\alpha \choose n} =
\left\{
   \begin{array}{ll}
 \displaystyle{\frac{\Gamma (n+\alpha+1)}{\Gamma(n+1)\Gamma(\alpha+1)},}
 & \mbox{if $n\in\N$, $\alpha\notin \Z_-$}, \\
 \displaystyle{\frac{(-1)^n \Gamma (-\alpha)}{\Gamma(n+1)\Gamma(-n-\alpha)},}
 & \mbox{if $n\in\N$, $n+\alpha\in \Z_-$}, \\
    0,
 & \mbox{if $n\in\N$, $\alpha \in \Z_-$, $n+\alpha\in \N_0$}, \\
   1,
 & \mbox{if $n=0$}, \\
   0,
 & \mbox{if $n\in\Z_-$}. \\
   \end{array}\right.
\label{def:combination}
\end{equation}
Note that for $n\in\N, \alpha\in\Z_-$ with $n+\alpha\le -1$,
$\displaystyle{
{n+\alpha \choose n}
=(-1)^n{-\alpha-1 \choose n}.}
$
By this definition, the equality
\begin{equation}
\left. \frac{1}{n!} \left( \frac{d}{dx}\right)^n x^{n+\alpha} \right|_{x=1}
= {n+\alpha \choose n}
\label{eqn:d/dx}
\end{equation}
holds for $n\in \N_0, \ \alpha\in\R$.
Then Laguerre polynomials can be expressed as
\begin{equation}
L_j^{\alpha}(x)
=\sum_{\ell=0}^j \frac{(-1)^\ell}{\ell !}
{j+\alpha \choose j-\ell}x^\ell
\label{eqn:Laguerre}
\end{equation}
for any $\alpha \in \R$.
Remark that 
applying (\ref{eqn:d/dx}) to the equation
\begin{equation}
\frac{1}{n!} \left( \frac{d}{dx} \right)^{n} x^{n+\alpha}
=\frac{1}{(n-1)!} \left( \frac{d}{dx} \right)^{n-1} x^{(n-1)+\alpha}
+ x \frac{1}{n!} \left( \frac{d}{dx} \right)^{n} x^{n+(\alpha-1)},
\nonumber
\end{equation}
with $x=1$ and putting $\beta=\alpha +n$, we have the identity
\begin{equation}
{\beta \choose n} = {\beta -1 \choose n}
+{\beta-1 \choose n-1},
\quad n\in\Z, \ \beta\in\R.
\label{formula:combination}
\end{equation}

We introduce the polynomials
\begin{eqnarray}
\label{eqn:F1}
F_j(x) &=& -\frac{d}{dx} L_{j+1}^{2\ma}(x), 
\quad j\in\N_0,
\\
\label{eqn:G1}
G_{j}(x) &=& \frac{d}{dx} \left\{
L_{j+1}^{2\ma}(x) - \frac{j+2\ma}{j} L_{j-1}^{2\ma}(x)\right\},
\quad j\in \N.
\end{eqnarray}
For $k\in\N_0, j=0,1,2,\dots,k$, let
\begin{eqnarray}
\alpha_{k, j}
&=& {k-j+\mb \choose k-j}, \qquad \text{ if $k$ is even},
\nonumber\\
\alpha_{k,j}
&=& \frac{k+2\ma}{k}{k-2-j+\mb \choose k-2-j}-{k-j+\mb \choose k-j},
\quad \text{ if $k$ is odd},
\label{def:alpha}
\end{eqnarray}
In Appendix \ref{chap:Proofs4142}, 
we will give the proof of the following
lemmas.


\begin{lem}
\label{thm:lemB2}
\quad 
For $\ell \in \N_0$
\begin{eqnarray}
&&F_{2\ell}(x)=\sum_{j=0}^{2\ell}\alpha_{2\ell,j}L^{\nu}_j(x),
\label{eqn:lemB2_1}
\\
&&G_{2\ell+1}(x)=\sum_{j=0}^{2\ell+1}\alpha_{2\ell +1,j}L^{\nu}_j(x).
\label{eqn:lemB2_2}
\end{eqnarray}

\end{lem}

\begin{lem}
\label{lem:skew0}
For $q, \ell\in \N_0$
\begin{eqnarray}
&& \langle F_{2q}, G_{2\ell+1} \rangle_{*}
= - \langle G_{2\ell+1}, F_{2q} \rangle_{*}
=r_{q}^{*} \delta_{q \ell}, \\
\label{eqn:sskew2}
&& \langle F_{2q}, F_{2\ell} \rangle_{*}=0,
\label{eqn:44_1}
\\
&&\langle G_{2q+1}, G_{2\ell+1} \rangle_{*}=0,
\label{eqn:44_2}
\end{eqnarray}
with
\begin{equation}
r_{q}^{*}\equiv \frac{4\Gamma(2q+2\ma+2)}{(2q+1)!}
=4\Gamma(2\ma+1){2q+2\ma+1 \choose 2q+1}.
\label{eqn:rstar}
\end{equation}
\end{lem}

Then if we define the monic polynomials in $x$
of degree $k$ for $k \in \N_0$ as
\begin{equation}
\label{def:Rk}
R_{k}(x) = k! \left( \frac{c_1}{ \chi_1} \right)^k
\sum_{j=0}^k \alpha_{k, j}
L^{\nu}_j \left( \frac{x}{c_1}\right) \chi_1^{j},
\end{equation}
Lemma \ref{lem:skew0} gives 
the following through the relation 
(\ref{eqn:fg_1}) and
the orthogonality of the Laguerre polynomials 
(\ref{eqn:orth}). 

\begin{lem}
\label{lem:skew1}
For $q, \ell \in \N_0$
\begin{eqnarray*}
\label{eq:skew1}
{\langle {R}_{2q}, {R}_{2\ell+1} \rangle }
=- {\langle {R}_{2\ell+1}, {R}_{2q} \rangle}
={r}_q \delta_{q \ell},
\\
\label{eq:skew2}
{\langle {R}_{2q}, {R}_{2\ell} \rangle}
= 0, \quad
{\langle {R}_{2q+1}, {R}_{2\ell+1} \rangle} = 0,
\end{eqnarray*}
where
\begin{equation}
r_q= 2^{-2\nu}T^{-\kappa}\left(\frac{t_1^2}{T}\right)^{4q+1}
\frac{(2q)!\Gamma(2q+2+2\ma)}{\Gamma(\nu+1)^2}.
\label{eqn:rq}
\end{equation}
\end{lem}

The choice of the polynomials
$F_{j}(x)$ and $G_{j}(x)$ in (\ref{eqn:F1}) and (\ref{eqn:G1}),
and their explicit expansions in terms of 
the Laguerre polynomials (Lemma \ref{thm:lemB2})
are crucial,
since they enable us to determine the
appropriate skew-orthogonal polynomials
(Lemma \ref{lem:skew1}).
As shown below, we are able to inverse the
skew-symmetric matrix $A_0$ given by
(\ref{eqn:A0B})
readily for arbitrary (even) $N$, by using
these skew-orthogonal polynomials.


\vskip 3mm

\subsection{Matrix inversion}

Let $b_{2k}=b_{2k+1}=r_{k}^{-1/2}, k \in \N_0$, and
determine the polynomials $\{M_{i}(x)\}_{0 \leq i \leq N-1}$
in (\ref{eqn:M}) as
\begin{equation*}
  M_{i}(x)=b_{i} R_{i}(x), \quad i=0,1, \cdots, N-1.
\label{eqn:M2}
\end{equation*}
Then by (\ref{eqn:A0B}), (\ref{eqn:A0B2}) and 
Lemma \ref{lem:skew1}, we have the equality
\begin{equation}
 \langle i | \hat{J}_{N} | j \rangle
 =(J_{N})_{ij}, \quad i, j =1,2, \cdots, N,
\label{eqn:A0C}
\end{equation}
where
$J_{N}=I_{N/2} \otimes J_{2}$.
It is interesting to compare this result
with (\ref{eqn:hatJ}).
Since $J_{N}^{2}=-I_{N}$, we can immediately 
obtain the inversion matrix appearing in 
(\ref{eqn:DSI1}) as
\begin{equation}
\langle i | (\hat{J}_{N})^{\bigtriangleup} |j \rangle
= - (J_{N})_{ij}, \quad
i,j =1,2, \cdots, N.
\label{eqn:inverse1}
\end{equation}
If we consider a semi-infinite matrix
$$
J \equiv \lim_{N \to \infty} J_{N}
=\Big( \langle i | \hat{J} | j \rangle \Big)_{i,j \in \N},
$$
its inverse matrix may be given by
$$
J^{-1}= \Big(\langle i | \hat{J}^{\bigtriangleup} | j \rangle 
\Big)_{i,j \in \N}
=-J.
$$
Using expansions (\ref{epn:Rmk}) and (\ref{epn:Phimk})
with Lemmas \ref{thm:lemB2}, \ref{lem:skew0} and \ref{thm:lemC2},
we can show
$$
\langle m,x|\hat{p}|n,y \rangle
= \langle m,x | \hat{p} \hat{J} |i \rangle
\langle i | \hat{J}^{\bigtriangleup} | j \rangle
\langle j |n, y \rangle
$$ 
and so
$$
\langle m,x|\hat{p}\hat{J}\hat{p}|n,y \rangle
= \langle m,x | \hat{p} \hat{J} |i \rangle
\langle i | \hat{J}^{\bigtriangleup} | j \rangle
\langle j |\hat{J} \hat{p} |n, y \rangle.
$$
Then the equations (\ref{eqn:SItilde1}) 
are written as
\begin{eqnarray}
\tS^{m,n}(x,y)
&=& \left\{
\begin{array}{ll}
\langle m, x | \hat{p} \hat{J} |i \rangle
\langle i | (\hat{J}_{N})^{\bigtriangleup} | j \rangle
\langle j |n, y \rangle,
& \mbox{if $m \geq n$,} \\
& \\
-\langle m, x | \hat{p} \hat{J} | i \rangle
\langle i | 
(\hat{J}^{\bigtriangleup} - (\hat{J}_{N})^{\bigtriangleup})
 | j \rangle \langle j | n, y \rangle,
& \mbox{if $m < n$,}
\end{array}
\right. \nonumber\\
\label{eqn:SItilde2}
\tI^{m,n}(x,y) 
&=& \langle m, x | \hat{p} \hat{J} | i \rangle
\langle i | (\hat{J}^{\bigtriangleup}
-(\hat{J}_{N})^{\bigtriangleup})| j \rangle
\langle j | \hat{J} \hat{p} | n, y \rangle.
\end{eqnarray}

Now we introduce the notations,
just following the previous papers 
for multi-matrix models \cite{NF99,FNH99,Nag01}, as
\begin{eqnarray}
R_{i}^{(m)}(x) &\equiv& \frac{1}{b_{i}} \langle m, x | i+1 \rangle
\nonumber\\
\label{def:Rim}
&=& 
\int_0^\infty dy \ R_{i}(y) \p({t_1},y |0)
\p(t_m-t_1,x| y), \\
\Phi_{i}^{(m)}(x) &\equiv& 
- \frac{1}{b_{i}} \langle m, x| \hat{p} \hat{J} | i+1 \rangle
\nonumber\\
\label{def:Phi}
&=& \int_{0}^{\infty} dy \
R^{(m)}_{i}(y) F^{(m)}(y,x),
\end{eqnarray}
for $i=0,1, \cdots, N-1, m=1,2, \cdots, M+1$, where
\begin{equation}
F^{(m)}(x,y)
= \int_{0}^{\infty} dw \int_{0}^{w} dz \ 
\left| \begin{array}{cc}
 \p(T-t_m, z|x) & \p(T-t_m, w|x)
 \cr
 \p(T-t_m, z|y) & \p(T-t_m,w|y )
\end{array}
\right|.
\label{def:Fmm}
\end{equation}
It should be noted that
$R_{i}^{(1)}(x)=R_{i}(x)\p (t_1,x|0),
 0 \leq i \leq N-1$,
and $F^{(1)}(x,y)=F(x,y)$,
where $R_{i}(x)$ and
$F(x,y)$ were defined by
(\ref{def:Rk}) and (\ref{def:F}), respectively.
Then we arrive at the following explicit expressions
for the elements of matrix kernel (\ref{eqn:matrixA})
of our Pfaffian processes,
\begin{eqnarray}
&&D^{m,n}(x,y)=D^{m,n}_N(x,y)
= \sum_{\ell=0}^{(N/2)-1} 
\frac{1}{r_{\ell}}\Big[R_{2\ell}^{(m)}(x) R_{2\ell+1}^{(n)}(y) 
- R_{2\ell+1}^{(m)}(x) R_{2\ell}^{(n)}(y)\Big],
\nonumber\\
&&\tI^{m,n}(x,y)=\tI^{m,n}_N(x,y)
= - \sum_{\ell=N/2}^{\infty} 
\frac{1}{r_{\ell}}\Big[\Phi_{2\ell}^{(m)}(x) \Phi_{2\ell+1}^{(n)}(y)
- \Phi_{2\ell+1}^{(m)}(x) \Phi_{2\ell}^{(n)}(y)\Big],
\nonumber\\
&&S^{m,n}(x,y)=S^{m,n}_N(x,y)
= \sum_{\ell=0}^{(N/2)-1} 
\frac{1}{r_{\ell}}\Big[\Phi_{2\ell}^{(m)}(x) R_{2\ell+1}^{(n)}(y) 
- \Phi_{2\ell+1}^{(m)}(x) R_{2\ell}^{(n)}(y)\Big],
\label{def:DISmn}
\end{eqnarray}
and
\begin{eqnarray}
&&\tS^{m,n}(x,y)=\tS^{m,n}_N(x,y)
=S^{m,n}(x,y)- \p(t_{n}-t_m,y|x){\bf 1}_{(m < n)}.
\label{def:tSmn}
\end{eqnarray}

\SSC{Asymptotic Behavior of Correlation Functions}
\label{chap:Asymptotic}

In this section, we give the proof of our main theorem
(Theorem \ref{thm:main1}), by estimating the
$N \to \infty$ asymptotic of matrix kernel (\ref{eqn:matrixA})
of Theorem \ref{thm:finite}.
Elementary calculation needed for the estimation
are summarized in Appendix \ref{chap:Elementary}.
Here $a_N \sim b_N, N \to\infty$ means 
$a_N/b_N \to 1, N\to\infty$.
We assume that $T=N$, $t_m=T+s_m, 1 \leq m \leq M+1$ with
$s_1 < s_2 < \cdots < s_{M} < s_{M+1}=0$.
We put
\begin{eqnarray}
&&L_j^{\nu}(x, -s_m)=L^{\nu}_{j}\left(x \right)\chi_m^{j}.
\quad
\hbox{ and } \quad
\hL_j^{\nu}(x, s_m)=\frac{\Gamma(j+1)}{\Gamma(j+1+\nu)}
L^{\nu}_{j}\left(x \right)\chi_m^{-j}.
\label{eqn:LLhat}
\end{eqnarray}

\subsection{Asymptotics of $R_{k}(x)$ and $R_{k}^{(m)}(x)$}

Let
\begin{equation*}
\hR_k(x) 
= \frac{1}{k!}\left( \frac{t_1^2}{T}\right)^{-k} R_k(x)
= \sum_{j=0}^k \alpha_{k,j}L_j^{\nu}\left(\frac{x}{c_1}, -s_1\right).
\end{equation*}
Since $c_1 \sim N = T$,
\begin{eqnarray}
\hR_{2\ell}(x)
&\sim& I(2\ell, \mb ),
\nonumber\\
\hR_{2\ell+1}(x)
&\sim& \frac{2\ma}{2\ell+1}I(2\ell-1,\mb)
-I(2\ell+1,\mb-1)-I(2\ell,\mb-1)
\nonumber\\
&\sim& \frac{\ma}{\ell}I(2\ell,\mb)-2I(2\ell,\mb-1),
\quad N \to \infty,
\label{eq:I_2l+1}
\end{eqnarray}
where 
\begin{equation*}
I(q,c) \equiv \sum_{j=0}^{q} { q-j+c \choose q-j}
L_j^{\nu}\left(\frac{x}{N}, -s_1 \right)
=\sum_{j=0}^{q} { j+c \choose j}
L_{q-j}^{\nu}\left(\frac{x}{N}, -s_1 \right)
\end{equation*}
for $q \in\N$ and $c\in\R$.
We set
$$
    2 \ell = N \theta,
$$
and examine the asymptotic behavior of $I(2\ell,c)$
as $N\to\infty$ with some $\theta \in (0, \infty)$.
When $c\in\Z_-$, 
$\displaystyle{{j+c \choose j}=(-1)^{j}
{-c-1 \choose j}}$. Then
from (\ref{sim:K}) in Lemma \ref{thm:lemD2}
with $j=2 \ell$ ({\it i.e.} $\eta=1$ in
(\ref{eqn:setting})), we can easily see
\begin{eqnarray*}
I(2\ell,c)
&=& \sum_{j=0}^{-c-1}{j+c \choose j}
L_{2\ell-j}^{\nu}\left(\frac{x}{N}, -s_1 \right)
\sim\frac{(N\theta)^{c+\nu+1}}{(\theta x)^{\nu}}
\Knu^{(-c-1)}(\theta,1,x,-s_1),
\quad N \to \infty.
\end{eqnarray*}
This result is generalized to the following lemma.

\begin{lem}
\label{thm:lem51}
\quad For any $c\in \R$,
$\theta \in (0, \infty)$, we have
\begin{eqnarray}
I(2\ell,c)
&\sim& \frac{(N\theta)^{c+\nu+1}}{(\theta x)^{\nu}}
\Knu^{(-c-1)}(\theta,1,x,-s_1), 
\quad N \to \infty.
\label{lemma:I(2ell,c)}
\end{eqnarray}
\end{lem}

\noindent{\it Proof. } \quad
\begin{eqnarray*}
&&I(2\ell,c)= \sum_{p=0}^{2\ell}{p+c \choose p}
L^{\nu}_{2\ell-p}\left(\frac{x}{N}, -s_1 \right)
\nonumber\\
&&=\sum_{p=0}^{2\ell}{p+c \choose p}
\Big\{ 
L^{\nu}_{2\ell-p}\left(\frac{x}{N}, -s_1 \right)
-L^{\nu}_{2\ell}\left(\frac{x}{N}, -s_1\right)
\Big\}
+\sum_{p=0}^{2\ell}{p+c \choose p}
L^{\nu}_{2\ell}\left(\frac{x}{N}, -s_1 \right)
\nonumber\\
&&=\sum_{p=0}^{2\ell}{p+c \choose p}
\sum_{k=0}^{p-1}\Big\{ 
L^{\nu}_{2\ell-k-1}\left(\frac{x}{N}, -s_1 \right)
-L^{\nu}_{2\ell-k}\left(\frac{x}{N}, -s_1 \right)
\Big\}
+\sum_{p=0}^{2\ell}{p+c \choose p}
L^{\nu}_{2\ell}\left(\frac{x}{N}, -s_1 \right)
\nonumber\\
&&=\sum_{p=0}^{2\ell}{p+c \choose p}
\sum_{k=0}^{p-1}\sum_{q=0}^1 (-1)^{q+1} {1 \choose q}
L^{\nu}_{2\ell-k-q}\left(\frac{x}{N}, -s_1 \right)
+\sum_{p=0}^{2\ell}{p+c \choose p}
L^{\nu}_{2\ell}\left(\frac{x}{N}, -s_1 \right).
\end{eqnarray*}
Repeating this procedure, we have
\begin{equation}
I(2\ell,c)= \sum_{k=0}^\infty (-1)^k a_k(2\ell,c)
\sum_{q=0}^k (-1)^q {k \choose q}
L_{2\ell-q}^\nu\left(\frac{x}{N}, -s_1 \right)
\label{eqn:I(2l,c)-1}
\end{equation}
with
\begin{eqnarray}
&&a_k(2\ell,c)=  \sum_{p=0}^{2\ell}{p+c \choose p}
\sum_{j_1=0}^{p-1}\sum_{j_2=0}^{j_1-1}\cdots
\sum_{j_k=0}^{j_{k-1}-1}1
=\sum_{p=0}^{2\ell}{p+c \choose p}
{p \choose k}.
\nonumber
\end{eqnarray}
Using (\ref{formula:combination}), 
we can rewrite $a_k(2\ell,c)$ as
\begin{eqnarray*}
&&a_k(2\ell,c)=\sum_{p=0}^{2\ell}
\Big\{ {p+c+1 \choose p} - {p+c \choose p-1} \Big\}{p \choose k}
\\
&&\quad= \sum_{p=0}^{2\ell}
\bigg[ {p+c+1 \choose p}{p \choose k} 
- {p+c \choose p-1}{p-1 \choose k} 
+{p+c \choose p-1}\Big\{ {p-1 \choose k} - {p \choose k} \Big\}
\bigg]
\\
&&\quad= {2\ell+c+1 \choose 2\ell}{2\ell \choose k} 
-a_{k-1}(2\ell-1,c+1).
\end{eqnarray*}
Using this equation recursively, we obtain 
\begin{equation}
a_k(2\ell,c)
= \sum_{r=0}^k (-1)^r 
{2\ell+c+1 \choose 2\ell-r}{2\ell-r \choose k-r}
= {2\ell+c+1 \choose 2\ell-k}{c+k \choose k}.
\label{eqn:a_k(2l,c)}
\end{equation}
Thus (\ref{eqn:I(2l,c)-1}) with (\ref{eqn:a_k(2l,c)}) gives
\begin{equation*}
I(2\ell,c)=\sum_{k=0}^\infty (-1)^k 
{2\ell+c+1 \choose 2\ell-k}{c+k \choose k}
\sum_{q=0}^k (-1)^q {k \choose q}
L_{2\ell-q}^\nu\left(\frac{x}{N}, -s_1 \right).
\end{equation*}
By simple calculation with the estimate (\ref{sim:choose}) 
and (\ref{sim:K}) of Lemma \ref{thm:lemD2}, 
we obtain (\ref{lemma:I(2ell,c)})
through the expression (\ref{def:K^c}).
\qed
\vskip 3mm

From above asymptotic of $I(2\ell,c)$ with equations
(\ref{eq:I_2l+1})
we have the following proposition.
\vskip 3mm

\begin{prop}
\label{thm:prop52} \quad
{\rm (1)} \ Suppose that $\ell\in\N$ and 
$2\ell \sim N\theta$, $N\to\infty$ for some $\theta\in (0,\infty)$.
Then
\begin{equation*}
\hR_{2\ell}(x)\sim 
\frac{(N\theta)^{2\ma+1}}{(\theta x)^{\nu}}
\Knu^{(-\mb-1)}(\theta,1,x,-s_1),
\quad N \to \infty.
\label{eqn:R-2l}
\end{equation*}

\noindent
{\rm (2)} \ Suppose that $\ell\in\N_0$ and 
$2\ell+1 \sim N\theta$, $N\to\infty$ for some $\theta\in (0,\infty)$.
Then
\begin{equation*}
\hR_{2\ell+1}(x)
\sim \frac{2(N\theta)^{\mb+\nu}}{(\theta x)^{\nu}}
\bigg[\ma\Knu^{(-\mb-1)}(\theta,1,x,-s_1)
-\Knu^{(-\mb)}(\theta,1,x,-s_1)\bigg],
\quad N \to \infty.
\label{eqn:R-2l+1}
\end{equation*}

\end{prop}

\vskip 3mm

We next examine asymptotic of $R_{k}^{(m)}(x)$.
From the definition (\ref{def:Rim}) and the expression
(\ref{eqn:expansionL})
\begin{eqnarray}
&&R^{(m)}_k(x) = c_1 \int_0^\infty d\eta \,
R_k(c_1 \eta) \p(t_1,c_1\eta |0)
\p(t_m-t_1,x |c_1 \eta)
\nonumber\\
&&\quad=
\frac{k!}{2^{\nu+1}\Gamma(\nu+1)}
\left( \frac{c_1}{\chi_1}\right)^k \left( \frac{1}{t_m}\right)^{\nu+1}
x^\ma\exp\left[ \left(-1+\frac{t_m}{2T}\right)\frac{x}{c_m}\right]
\sum_{j=0}^k \alpha_{k,j} L_j^{\nu}\left(\frac{x}{c_m}, -s_m \right).
\nonumber\\
\label{epn:Rmk}
\end{eqnarray}
We put 
\begin{equation}
\label{def:hRkm}
\hR_k^{(m)}(x) 
= \frac{2^{\nu}T^{\nu}\Gamma(\nu+1)}{\Gamma(k+1+2\ma)}
\left( \frac{\chi_1}{c_1}\right)^k 
R_k^{(m)}(x).
\end{equation}
If we set $k \sim N \theta$ as $N \to \infty$,
(\ref{sim:Gamma/Gamma}) in Appendix \ref{chap:Elementary} gives
$$
\frac{k! T^{\nu}}{\Gamma(k+1+2\ma)}\left( \frac{1}{t_m}\right)^{\nu+1}
\exp\left[ \left(-1+\frac{t_m}{2T}\right)\frac{x}{c_m}\right]
\sim N^{-(2\ma+1)}\theta^{-2\ma}, \quad N\to\infty,
$$
then we obtain the following from Proposition \ref{thm:prop52}.

\vskip 3mm

\begin{prop}
\label{thm:prop53}
\quad
{\rm (1)} \ Suppose that $\ell\in\N$ and 
$2\ell \sim N\theta$, $N\to\infty$ for some $\theta\in (0,\infty)$.
Then 
\begin{equation*}
\hR^{(m)}_{2\ell}(x)\sim 
\frac{\theta^{1-\nu} x^{-\kappa/2}}{2}
\Knu^{(-\mb-1)}(\theta,1,x,-s_m),
\quad N \to \infty.
\label{eqn:Rm-2l}
\end{equation*}

\noindent
{\rm (2)} \ Suppose that $\ell\in\N_0$ and 
$2\ell+1 \sim N\theta$, $N\to\infty$ for some $\theta\in (0,\infty)$.
Then
\begin{equation*}
\hR^{(m)}_{2\ell+1}(x)\sim 
\frac{\theta^{-\nu} x^{-\kappa/2}}{N}
\bigg[ 
\ma\Knu^{(-\mb-1)}(\theta,1,x,-s_m)-\Knu^{(-\mb)}
(\theta,1,x,-s_m)\bigg],
\quad N \to \infty.
\label{eqn:Rm-Rm-2l+1}
\end{equation*}
\end{prop}


\subsection{Asymptotics of $\Phi_{k}^{(m)}(x)$}

Using (\ref{eqn:gmM+1}), 
(\ref{def:Fmm}) is rewritten as
\begin{eqnarray*}
F^{(m)}(y,x)
&=& \left(\frac{1}{T}\right)^{\kappa}
\left(\frac{t_m}{c_m}\right)^{2(\nu+1)}
(xy)^{\kappa/2}
\exp\left\{\left(-\frac{t_m}{T}\right)\frac{x+y}{2c_m}\right\}
\nonumber\\
&&\qquad \times
\sum_{p=0}^{\infty}\sum_{j=0}^{\infty}
\langle L^{\nu}_p, L^{\nu}_j \rangle_*
\hL^{\nu}_p \left(\frac{y}{c_m}, s_m \right)
\hL^{\nu}_j \left(\frac{x}{c_m}, s_m \right).
\label{eqn:Fmm}
\end{eqnarray*}
Hence from (\ref{def:Phi}) and (\ref{epn:Rmk}), we have
\begin{eqnarray}
\Phi_k^{(m)}(x) 
&=&c_m \int_0^{\infty} d\eta \, R_k^{(m)}(c_m\eta)
F^{(m)}(c_m \eta, x)
\nonumber\\
&=& 
\frac{k!}{2^{\nu+1}\Gamma(\nu+1)}\left(\frac{c_1}{\chi_1}\right)^k
\left( \frac{1}{c_m} \right)^{\nu+1}
\frac{t_m^{\nu+1}}{T^{\kappa}}
x^{\kappa/2}
\exp\left\{ -\frac{t_m x}{2Tc_m}\right\}
\nonumber\\
&&\qquad \times \sum_{j=0}^\infty 
\langle \sum_{p=0}^k \alpha_{k,p}L_p^{\nu},
L_j^{\nu} \rangle_*
\hL_j^{\nu}\left(\frac{x}{c_m}, s_m \right),
\label{epn:Phimk}
\end{eqnarray}
where we have used the orthogonal relation (\ref{eqn:orth}) 
of Laguerre polynomials.
Put
\begin{equation}
\hPhi_k^{(m)}(x) 
= \frac{2^{\nu}T^{-\mb}\Gamma(\nu+1)}{k!}
\left(\frac{\chi_1}{c_1}\right)^k\Phi_k^{(m)}(x).
\label{eqn:phihat}
\end{equation}
Then we have the following proposition.

\vskip 3mm
\begin{prop}
\label{thm:prop54}
\quad {\rm (1)} \ Suppose that $\ell\in\N$ and 
$2\ell \sim N\theta$, $N\to\infty$ for some $\theta\in (0,\infty)$.
Then
\begin{equation}
\hPhi_{2\ell}^{(m)}(x) 
\sim -\theta^{\nu}x^{\kappa/2}
\int_1^\infty d\xi \ \xi^\ma \hKnu^{(\mb+1)}(\theta,\xi,x,s_m),
\quad \N \to \infty.
\label{Phi_2l}
\end{equation}

\noindent
{\rm (2)} \ Suppose that $\ell\in\N_0$ and 
$2\ell+1 \sim N\theta$, $N\to\infty$ for some $\theta\in (0,\infty)$.
Then
\begin{equation}
\hPhi_{2\ell+1}^{(m)}(x) \sim 
-\frac{2\theta^{-1+\nu}x^{\kappa/2}}{N}
\hKnu^{(\mb+1)}(\theta,1,x,s_m),
\quad N \to \infty.
\label{eqn:Phi2l+1}
\end{equation}

\end{prop}

\noindent{\it Proof.}
\begin{eqnarray*}
\hPhi_k^{(m)}(x) 
&\sim& \frac{T^{-\nu}x^{\kappa/2}}{2}
\sum_{j=0}^{\infty} 
\langle \sum_{p=0}^k \alpha_{k,p}L_p^{\nu},L_j^{\nu}  \rangle_*
\hL_j^{\nu}\left( \frac{x}{N}, s_m \right),
\quad N\to\infty,
\\
&=& \frac{T^{-\nu}x^{\kappa/2}}{2}
\sum_{j=0}^{\infty} 
\langle  Q_k, \sum_{q=0}^j \beta_{j,q}Q_q \rangle_*
\hL_j^{\nu}\left( \frac{x}{N}, s_m \right).
\end{eqnarray*}
Here we have used the notation (\ref{def:H})
and introduced $\vbeta = (\beta_{j, k})$, which is the inverse
of the matrix $\valpha =(\alpha_{k,j})$ given by
Lemma \ref{thm:lemC2} in Appendix \ref{chap:beta}.
By the skew orthogonality of 
$\{Q_{k}\}$ given by Lemma \ref{lem:skew0}, we have
\begin{eqnarray}
&&\hPhi_{2\ell}^{(m)}(x) 
\sim  \frac{T^{-\nu}x^{\kappa/2}r_{\ell}^*}{2}
\sum_{j=0}^{\infty} \beta_{j,2\ell+1}
\hL_j^{\nu}\left( \frac{x}{N}, s_m \right),
\quad N\to\infty,
\label{hphi_2l}
\\
&&\hPhi_{2\ell+1}^{(m)}(x) 
\sim-\frac{T^{-\nu}x^{\kappa/2}r_{\ell}^*}{2}
\sum_{j=0}^{\infty} \beta_{j,2\ell}
\hL_j^{\nu}\left( \frac{x}{N}, s_m \right),
\quad N\to\infty.
\label{hphi_2l+1}
\end{eqnarray}


By (\ref{eqn:rstar}) and (\ref{def:beta_even}),
(\ref{hphi_2l+1}) gives
\begin{eqnarray}
&&\hPhi_{2\ell+1}^{(m)}(x) 
\sim -2\Gamma(2\ma+1)T^{-\nu}x^{\kappa/2}
{2\ell+1+2\ma \choose 2\ell +1}
\nonumber\\
&&\quad\times
\sum_{j=2\ell}^{\infty}
{j-2\ell-\mb-2 \choose j-2\ell}
\hL^{\nu}_{j}\left(\frac{x}{N}, s_m \right),
\quad N\to\infty.
\label{eqn:hphi_2l+1_2}
\end{eqnarray}
From (\ref{formula:combi_2}) we have
\begin{eqnarray}
&&\sum_{j=2\ell}^{\infty}
{j-2\ell-\mb-2 \choose j-2\ell}
\hL^{\nu}_{j}\left(\frac{x}{N}, s_m \right)
= \sum_{r=0}^{\infty}
{r-\mb-2 \choose r}
\hL^{\nu}_{2\ell+r}\left(\frac{x}{N}, s_m \right)
\nonumber\\
&&\qquad= \sum_{r=0}^{\infty}
{r-\mb-2+\alpha \choose r}
\sum_{p=0}^\alpha (-1)^p {\alpha \choose p}
\hL^{\nu}_{2\ell+r+p}\left(\frac{x}{N}, s_m \right)
\nonumber\\
&&\qquad\sim
(2\ell)^{-\alpha}
\sum_{r=0}^{\infty}
{r-\mb-2+\alpha \choose r}
\hKnu^{(\alpha)}(\theta,\eta,x,s_m),
\quad N\to\infty,
\label{eqn:sum_tL=sum_K}
\end{eqnarray}
where (\ref{sim:hK}) of Lemma \ref{thm:lemD2} was applied.
Setting $j=2 \ell \eta = N \theta \eta$ and using (\ref{sim:choose})
in Appendix \ref{chap:Elementary},
we conclude that 
\begin{eqnarray}
\hPhi_{2\ell+1}^{(m)}(x) 
&\sim&
-2 T^{-\nu}(2\ell)^{2\ma-\alpha}x^{\kappa/2}
\sum_{j=2\ell}^{\infty}
\frac{(j-2\ell+1)^{-\mb-2+\alpha}}{\Gamma(-\mb-1+\alpha)}
\hKnu^{(\alpha)}(\theta,\eta,x,s_m)
\nonumber\\
&\sim&
-\frac{2\theta^{-1+\nu}x^{\kappa/2}}{N\Gamma (-\mb-1+\alpha)} 
\int_1^{\infty} d\eta \
\frac{\hKnu^{(\alpha)}(\theta,\eta,x,s_m)}{(\eta-1)^{\mb+2-\alpha}},
\quad N\to\infty.
\end{eqnarray}
Through the expression (\ref{def:hK^c}),
we obtain (\ref{eqn:Phi2l+1}).


By (\ref{def:beta_odd}) of Lemma \ref{thm:lemC2}
with (\ref{eqn:bmn})
\begin{eqnarray*}
&&\sum_{j=2\ell+1}^{\infty}
\beta_{j,2\ell+1}
\hL^{\nu}_{j}\left(\frac{x}{N}, s_m \right)
\nonumber\\
&&\qquad= \frac{-1}{b(1,2\ell +1)}
\sum_{j=2\ell+1}^{\infty} 
\hL^{\nu}_j\left(\frac{x}{N}, s_m \right)
\sum_{r=\ell+1}^{[(j+1)/2]} b(1,2r-1){j-2r-\mb-1 \choose j-2r+1}
\nonumber\\
&&\qquad= \frac{-1}{b(1,2\ell +1)}S(\ell),
\end{eqnarray*}
where
\begin{equation*}
S(\ell)=\sum_{r=\ell+1}^{\infty} b(1,2r-1)
\sum_{j=2r-1}^{\infty}
\hL^{\nu}_j\left(\frac{x}{N}, s_m \right){j-2r-\mb-1 \choose j-2r+1}.
\end{equation*}
By this equation with the estimate (\ref{sim:b1l})
for $b(1, 2r-1)$ and (\ref{eqn:rstar}),
(\ref{hphi_2l}) becomes
\begin{eqnarray}
\hPhi_{2\ell}^{(m)}(x) &\sim&
-2T^{-\nu}(2\ell+2)^{-\ma} \Gamma(2\ma+1)
{2\ell+1+2\ma \choose 2\ell+1}x^{\kappa/2}S(\ell)
\nonumber\\
&\sim&
-2T^{-\nu} (2\ell +2)^\ma x^{\kappa/2}S(\ell),
\quad N\to\infty.
\label{eqn:hphi_2l_2}
\end{eqnarray}
From (\ref{eqn:sum_tL=sum_K}) with (\ref{sim:b1l})
\begin{eqnarray*}
S(\ell) &\sim&
\sum_{r=\ell+1}^{\infty} (2r)^\ma
\left(\frac{1}{2\ell}\right)^{\alpha}
\sum_{j=2r-1}^{\infty}
\frac{(j-2r+1)^{-\mb-2+\alpha}}{\Gamma(-\mb-1+\alpha)}
\hKnu^{(\alpha)}(\theta,\eta,x,s_m)
\nonumber\\
&\sim&
\frac{(2\ell)^{\kappa/2}}{2\Gamma(-\mb-1+\alpha)}
\int_1^{\infty} d\xi \ \xi^{\ma}
\int_{\xi}^{\infty} d\eta \ 
\frac{\hKnu^{(\alpha)}(\theta,\eta,x,s_m)}{(\eta-\xi)^{\mb+2-\alpha}},
\quad N\to\infty.
\end{eqnarray*}
Thus
\begin{eqnarray}
\hPhi_{2\ell}^{(m)}(x) 
&\sim& -\frac{\theta^{\nu}x^{\kappa/2}}{\Gamma(-\mb-1+\alpha)}
\int_1^{\infty} d\xi \ \xi^{\ma}
\int_{\xi}^{\infty} d\eta \ 
\frac{\hKnu^{(\alpha)}(\theta,\eta,x,s_m)}{(\eta-\xi)^{\mb+2-\alpha}}
\nonumber\\
&=& -\theta^{\nu}x^{\kappa/2}
\int_1^{\infty} d\xi \ \xi^{\ma}
\hKnu^{(\mb+1)}(\theta,\xi,x,s_m),
\quad N\to\infty,
\nonumber
\end{eqnarray}
where we used the expression (\ref{def:hK^c}).
This completes the proof of Proposition \ref{thm:prop54}.
\qed

\subsection{Asymptotics of $D^{m,n}(x,y)$, $\tI^{m,n}(x,y)$,
$S^{m,n}(x,y)$ and $\tS^{m,n}(x,y)$}

From the expressions (\ref{def:DISmn}) with (\ref{eqn:rq})
and the definitions (\ref{def:hRkm}) 
and (\ref{eqn:phihat}) we have 
\begin{eqnarray*}
&&D^{m,n}_N(x,y)\sim \sum_{\ell=0}^{(N/2)-1}
\left( \frac{2\ell}{N}\right)^{2\ma}
\Big[ \hR_{2\ell}^{(m)}(x)\hR_{2\ell+1}^{(n)}(y) 
- \hR_{2\ell+1}^{(m)}(x)\hR_{2\ell}^{(n)}(y) 
\Big],
\\
&&\tI^{m,n}_{N}(x,y)
\sim -\sum_{\ell=(N/2)}^{\infty}
\left( \frac{2\ell}{N}\right)^{-2\ma}
\Big[ \hPhi_{2\ell}^{(m)}(x) \hPhi_{2\ell+1}^{(n)}(y) 
- \hPhi_{2\ell+1}^{(m)}(x) \hPhi_{2\ell}^{(n)}(y) 
\Big],
\\
&&S^{m,n}_{N}(x,y)\sim \sum_{\ell=0}^{(N/2)-1}
\Big[ \hPhi_{2\ell}^{(m)}(x)\hR_{2\ell+1}^{(n)}(y) 
- \hPhi_{2\ell+1}^{(m)}(x)\hR_{2\ell}^{(n)}(y) 
\Big], \quad N \to \infty.
\end{eqnarray*}
From Propositions \ref{thm:prop53} and \ref{thm:prop54} 
we obtain the following asymptotics: 
\begin{eqnarray*}
D^{m,n}_N(x,y) &\sim& \cD (s_m, x; s_n,y),
\\
\tI^{m,n}_N(x,y) &\sim& \cI (s_m,x; s_n,y),
\\
S^{m,n}_N(x,y) &\sim& \cS (s_m,x; s_n, y),
\quad N\to\infty,
\end{eqnarray*}
where $\cD, \cI, \cS$ are defined by
(\ref{eqn:mainDIS}).

Next we study the asymptotic behavior of 
$\p (t_n-t_m, y|x)$.
From (\ref{eqn:expansionL}) we have
\begin{eqnarray*}
&&\p(t_n-t_m, y|x) 
\\
&&\qquad=\left(\frac{t_m}{t_n}\right)^{\nu+1}
c_m^{-\ma-1}\left(\frac{x}{c_m}\right)^{\kappa/2}y^\ma
\exp \left[-\frac{t_mx}{2Tc_m}\right]
\exp\left[ \left( -2+\frac{t_n}{T}\right) \frac{y}{2c_n}\right]
\\
&&\qquad\times\sum_{j=0}^{\infty}
\hL_j^\nu\left(\frac{x}{c_m}, s_m \right)
L_j^\nu\left(\frac{y}{c_n}, -s_n \right).
\end{eqnarray*}
Then by simple calculation with 
Lemma \ref{thm:lemD2} with $\alpha=0$, we have
\begin{eqnarray*}
\p(t_n-t_m, y|x) &\sim&
\left(\frac{y}{x}\right)^{\mb/2}
\frac{1}{N}\sum_{j=0}^{\infty}
\exp\left[ 2(s_m-s_n)\theta\eta\right]
J_\nu(2\sqrt{\theta\eta x})J_\nu(2\sqrt{\theta\eta y}),
\\
&\sim& \left(\frac{y}{x}\right)^{\mb/2}\cG(s_m,x;s_n, y),
\quad N\to\infty,
\label{sim:cp}
\end{eqnarray*}
where $\cG$ is defined by (\ref{eqn:mainG}), 
and then
$
\tS^{m,n}_N(x,y) \sim \tcS (s_m,x; s_n, y), \, N\to\infty.
$
Then, the proof of Theorem \ref{thm:main1} is completed.

\clearpage
\appendix

{\LARGE{\bf Appendices}}
\SSC{Proof of (\ref{eqn:Psi2})}\label{chap:Proof(3.16)}

We assume that the number of particles $N$ is even. 
Consider the multiple integral
\begin{eqnarray}
\kZ[\chi] &=& \left(\frac{1}{N !}\right)^{M+1}
\int_{\Rp^{N(M+1)}} \prod_{m=1}^{M+1} d \x^{(m)} \,
\det_{1 \leq i, j \leq N} \left[
M_{i-1}(x_{j}^{(1)}) \p(t_{1}, x_{j}^{(1)}|0)
(1+ \chi_{1}(x_{j}^{(1)})) \right] \nonumber\\
&& \times
\prod_{m=1}^{M} \det_{1 \leq i, j \leq N}
\left[ \p(t_{m+1}-t_{m}, x_{j}^{(m+1)}| x_{i}^{(m)})
(1+\chi_{m+1}(x_{j}^{(m+1)})) \right]
{\rm sgn}\left(h_{N}(\x^{(M+1)}) \right).
\nonumber
\end{eqnarray}
By the definition (\ref{eqn:Psi0}) with (\ref{eqn:expect})
and (\ref{eqn:gNT2}), and by the equality (\ref{eqn:hN2}),
we have
\begin{equation}
\label{eqn:Psi1}
\kPsi ( {\bf f}; \vtheta )=
\frac{\kZ[\chi]}{\kZ[0]},
\end{equation}
where
$\kZ[0]$ is obtained from $\kZ[\chi]$ by setting
$\chi_{m}(x) \equiv 0$ for all $m=1,2, \cdots, M+1$.

By repeated applications of the Heine identity
$$
\int_{\RW} d \x \,
\det_{1 \leq i, j \leq N} \Bigg[\phi_{i}(x_{j})\Bigg]
\det_{1 \leq i, j \leq N} \Bigg[ \bar{\phi}_{i}(x_{j}) \Bigg]
= \det_{1 \leq i, j \leq N} \left[
\int_{\Rp} dx \, \phi_{i}(x) \bar{\phi}_{j}(x) \right]
$$
for square integrable continuous functions
$\phi_{i}, \bar{\phi}_{i}, 1 \leq i \leq N$,
we have
\begin{eqnarray}
\kZ[\chi] &=& \int_{\RW} d \y \, 
\det_{1 \leq i, j \leq N} \left[
\int_{\Rp^{M+1}} \prod_{m=1}^{M+1} dx^{(m)} \,
\Bigg\{ M_{i-1}(x^{(1)}) \p(t_{1}, x^{(1)}|0)
(1+ \chi_{1}(x^{(1)})) \Bigg\} \right. \nonumber\\
&& \quad \left. \times \prod_{m=1}^{M} \Bigg\{
\p(t_{m+1}-t_{m}, x^{(m+1)}| x^{(m)} ) 
(1+\chi_{m+1}(x^{(m+1)}))
\Bigg\}
\delta (y_{j}-x^{(M+1)}) \right].
\nonumber
\end{eqnarray}
Using the notations in Section \ref{section32}, 
it is expressed as
\begin{eqnarray}
\kZ[\chi] 
&=& \int_{\RW} d \y \, 
\det_{1 \leq i, j \leq N} \left[
\langle i | \left(
1+ \frac{1}{1-\hat{\chi} \hat{p}_{+}}
\hat{\chi} \hat{p} \right)_{N} | M+1, y_j \rangle
\right]
\nonumber\\
&=& \int_{\RW} d \y \, 
\det_{1 \leq i, j \leq N} \left[
\langle M+1, y_{i} | \left(
1+ \hat{p} \hat{\chi} 
\frac{1}{1-\hat{p}_{-} \hat{\chi}} \right)_{N}
 | j \rangle
\right], \nonumber
\end{eqnarray}
since $\langle m, x | 
(\hat{p}_{+})^{k} | n, y \rangle
=\langle n, y | (\hat{p}_{-})^{k} |
m, x \rangle \equiv 0$
for $k > n-m \geq 0$.
Here we have used the Chapman-Kolmogorov equation,
$
\int_{\Rp} dy \, \p(t-s, y|x) \p(u-t, z|y)=
\p(u-s, z|x)
$,
$0 \leq s \leq t \leq u \leq T, \, x, y \in \Rp$.
Next we use the formula of de Bruijin \cite{Bruijn55}
$$
\int_{\RW} d\y \,
\det_{1 \leq i, j \leq N} \Bigg[
\phi_{i}(y_{j}) \Bigg]
= \Pf_{1 \leq i, j \leq N} \left[
\int_{\Rp} dy \int_{\Rp} d \tilde{y} \,
{\rm sgn}(\tilde{y}-y) \phi_{i}(y) \phi_{j}(\tilde{y})
\right]
$$
for integrable continuous functions
$\phi_{i}, 1 \leq i \leq N$,
in which the Pfaffian is defined by (\ref{def:pfaffian}).
Since $(\Pf(A))^2=\det A$ 
for any even-dimensional skew-symmetric matrix $A$, we have
\begin{eqnarray}
\Bigg( \kZ[\chi] \Bigg)^2
&=& \det_{1 \leq i, j \leq N} \left[
\langle i | 
\left( 1+ \frac{1}{1-\hat{\chi} p_{+}}
\hat{\chi} \hat{p} \right)_{N} \hat{J}
\left( 1+ \hat{p} \hat{\chi}
\frac{1}{1-\hat{p}_{-} \hat{\chi}} \right)_{N}
| j \rangle
\right] \nonumber\\
&=& \det_{1 \leq i, j \leq N} \Bigg[
(A_{0})_{ij}+(A_{1})_{ij}+(A_{2})_{ij}+(A_{3})_{ij} \Bigg]
\nonumber
\end{eqnarray}
with
\begin{eqnarray}
\label{eqn:A0}
(A_{0})_{ij} 
&=& \langle i | \hat{J}_{N} |j \rangle, 
\nonumber\\
(A_{1})_{ij} &=&
\langle i | \Big(
\frac{1}{1-\hat{\chi} \hat{p}_{+}}
\hat{\chi} \hat{p} \hat{J}\Big)_{N}| j \rangle
= \langle i | \Big( \hat{\chi}
\frac{1}{1-\hat{p}_{+} \hat{\chi}}
\hat{p} \hat{J}\Big)_{N}| j \rangle, 
\nonumber\\
(A_{2})_{ij} &=& 
\langle i | \Big( \hat{J} \hat{p}
\hat{\chi} \frac{1}{1-\hat{p}_{-}\hat{\chi}}
\Big)_{N}|
j \rangle, 
\nonumber\\
(A_{3})_{ij} &=& \langle i | \Big( \hat{\chi}
\frac{1}{1-\hat{p}_{+}\hat{\chi}}
 \hat{p}\hat{J} \hat{p} \hat{\chi}
 \frac{1}{1-\hat{p}_{-}\hat{\chi}} \Big)_{N} | j \rangle.
 \nonumber
\end{eqnarray}

Since $\Big( Z_{N,T}[0] \Big)^2
=\det_{1 \leq i, j \leq N} \Big[ (A_{0})_{ij} \Big]$, 
(\ref{eqn:Psi1}) gives
\begin{equation}
\label{eqn:det1}
\Bigg\{ \Psi_{N,T} ( {\bf f}; \vtheta ) 
\Bigg\}^2
= \det_{1 \leq i, j \leq N} \Bigg[
\delta_{i j}+(A_{0}^{-1}A_{1})_{ij}+(A_{0}^{-1}A_{2})_{ij}
+(A_{0}^{-1}A_{3})_{ij} \Bigg].
\end{equation}
By our notation (\ref{eqn:Delta1}),
$
(A_{0}^{-1})_{ij}=
\langle i |
(\hat{J}_{N})^{\bigtriangleup} | j \rangle,
$
and it is easy to confirm that (\ref{eqn:det1}) 
is written in the form
\begin{eqnarray}
\Bigg\{ \kPsi( {\bf f}; \vtheta ) 
\Bigg\}^2
&=& \det_{1 \leq i, j \leq N}
\Bigg[ \delta_{i j}
+ \langle i | {\bf B} | m,x \rangle 
\langle m,x | {\bf C} |j \rangle \Bigg], 
\label{eqn:det1b}
\end{eqnarray}
where we have introduced ${\bf B}$
and ${\bf C}$
as the following two-dimensional row and
column vector-valued
operators,
\begin{eqnarray}
 {\bf B} 
&=&  \left( \begin{array}{lll}
(\hat{J}_{N})^{\bigtriangleup} \circ \hat{\chi}
\frac{1}{1-\hat{p}_{+} \hat{\chi}} 
 & \; &
-(\hat{J}_{N})^{\bigtriangleup}
\circ \Big(1+\hat{\chi} \frac{1}{1-\hat{p}_{+} \hat{\chi}}\hat{p}\Big)
\hat{J} \hat{p} \hat{\chi} \\
\end{array}
\right), \nonumber\\
{\bf C} 
&=& \left( 
\begin{array}{c}
\hat{p} \hat{J}
\cr\cr
- \frac{1}{1-\hat{p}_{-} \hat{\chi}}
\end{array}
 \right). \nonumber
\end{eqnarray}
The determinant (\ref{eqn:det1b}) is equivalent with
the Fredholm determinant,
$$
{\rm Det} \langle m,x | I_{2}+{\bf C}\circ {\bf B}  |n, y \rangle. 
$$
Introducing matrix-valued operators,
$$
{\bf K_+}
= \left(
\begin{array}{ccc}
 1 - \hat{p}_{+} \hat{\chi}& & 0 
\cr & & \cr
0 && 1
\end{array}
 \right),
\quad
{\bf K_-}
= \left(
\begin{array}{ccc}
 1 & & 0 
\cr & & \cr
0 & & 1- \hat{p}_{-} \hat{\chi} 
\end{array}
 \right),
\quad
\widehat{\bf K}
= \left(
\begin{array}{ccc}
 1 & & -\hat{p}\hat{J}\hat{p}\hat{\chi}
\cr & & \cr
0 & & 1
\end{array}
 \right),
$$
we have
\begin{eqnarray}
I_{2}+{\bf C}\circ {\bf B}
&=& {\bf K_-}^{-1} 
\left[ {\bf K_-}{\bf K_+}+ 
\left(
\begin{array}{ccc}
 \hat{p} \hat{J} {\cal J}^{\bigtriangleup}_N
&\;&
\hat{p} \hat{J} \Big(1-{\cal J}^{\bigtriangleup}_N \hat{J} \Big)
 \hat{p} \cr
&\;& \cr
-{\cal J}^{\bigtriangleup}_N
&\;& 
{\cal J}^{\bigtriangleup}_N \hat{J} \hat{p} 
\end{array}  \right)\hat{\chi}
\right] {\bf K_+}^{-1}\widehat{\bf K}
\nonumber\\
&=& {\bf K_-}^{-1} 
\left[ I_2 + 
\left(
\begin{array}{ccc}
 \hat{p} \hat{J} 
{\cal J}^{\bigtriangleup}_N -\hat{p}_{+}
&\;&
\hat{p} \hat{J}\hat{p} 
-\hat{p} \hat{J}{\cal J}^{\bigtriangleup}_N \hat{J} \hat{p} 
\cr &\;& \cr
-{\cal J}^{\bigtriangleup}_N
&\;& 
{\cal J}^{\bigtriangleup}_N\hat{J} \hat{p} -\hat{p}_{-}
\end{array}  \right)\hat{\chi}
\right] {\bf K_+}^{-1}\widehat{\bf K},
\nonumber
\end{eqnarray}
where 
${\cal J}^{\bigtriangleup}_N=\circ(\hat{J}_{N})^{\bigtriangleup}\circ$.
From the orthogonality (\ref{eqn:orthZ}) and the definitions
(\ref{eqn:p+-}) of the operators
$\hat{p}_{+}$ and $\hat{p}_{-}$,
we have the fact that
$$
{\rm Det} \langle m, x | {\bf K_+}| n, y \rangle
={\rm Det} \langle m, x | {\bf K_-}| n, y \rangle 
={\rm Det} \langle m, x | \widehat{\bf K}| n, y \rangle
=1.
$$
Then (\ref{eqn:Psi2}) is derived.

\SSC{Proofs of Lemmas \ref{thm:lemB2} and \ref{lem:skew0}}
\label{chap:Proofs4142}

We use the following 
orthogonal relations
and formulae on Laguerre polynomials, which hold
for $\alpha, \beta >-1$;
\begin{eqnarray}
\label{eqn:orth}
&& \int_0^{\infty}L_j^\alpha (x)L_k^\alpha (x)x^\alpha e^{-x}dx
=\frac{\Gamma(\alpha +j+1)}{\Gamma(j+1)}\delta_{j k},
\quad j, k \in \N_{0},\\ 
\label{eqn:form1}
&& x \frac{d}{dx} L_{j}^{\alpha}(x) = j L_{j}^{\alpha}(x)
- (j+\alpha) L_{j-1}^{\alpha}(x), 
\quad j\in\N,
\\
\label{eqn:form2}
&&  L_{j}^{\alpha}(x)= - \frac{d}{dx} L_{j+1}^{\alpha}(x)
+ \frac{d}{dx} L_{j}^{\alpha}(x),
\quad j\in\N_0,
\\
&& L_{j}^{\beta}(x)
=\sum_{k=0}^{j}{j-k+\beta-\alpha-1 \choose j-k}L_{k}^{\alpha}(x),
\quad j\in\N_0.
\label{eqn:form3}
\end{eqnarray}

\vskip 0.3cm
\noindent{\bf Remark 6.} \quad
The identities (\ref{eqn:form1}) and (\ref{eqn:form2}) 
are given as Eqs. (6.2.6) and (6.2.7) in \cite{AAR99}. 
The relation (\ref{eqn:form3}) is proved in [1] as
(6.2.37) only when $\beta \geq \alpha >-1$.
The identity
(see (\ref{def:combination}) and \cite{Rio79})
\begin{equation*}
\sum_{\ell=0}^k {\ell-\alpha-1 \choose \ell}
{k-\ell+\alpha-1 \choose k-\ell}
= \delta_{k 0},
\end{equation*}
can be used to invert the relation
(\ref{eqn:form3}) to the form
$$
L_{\ell}^{\alpha}(x)
=\sum_{j=0}^{\ell}{\ell-j+\alpha-\beta-1 \choose \ell-j}L_{j}^{\beta}(x),
\quad \ell\in\N_0.
$$
Therefore, the validity of (\ref{eqn:form3})
for $\beta \geq \alpha > -1$ implies
that for $\alpha > \beta > -1$.

\subsection{Proof of Lemma \ref{thm:lemB2}}

In this subsection we prove Lemma \ref{thm:lemB2},
which gives the expansion formulae
of $F_k(x)$ and $G_k(x)$ in terms of $\{L_{j}^{\nu}(x)\}$.
Taking the summation of the equalities (\ref{eqn:form2})
from $0$ to $k\in\N_0$ gives
\begin{equation}
\sum_{n=0}^{k} L_{n}^{2\ma}(x)
= - \frac{d}{dx} L_{k+1}^{2\ma}(x)+ \frac{d}{dx} L_{0}^{2\ma}(x)
=F_{k}(x).
\label{eqn:F2}
\end{equation}
From (\ref{eqn:form3}), (\ref{eqn:F2}) and (\ref{formula:combination}),
we have
\begin{eqnarray}
F_{k}(x) 
&=& L_{k}^{2\ma}(x) + F_{k-1}(x) 
\nonumber\\
&=& \sum_{j=0}^{k} {k-j+\mb-1 \choose k-j}L_{j}^{\nu}(x)+ F_{k-1}(x) 
\nonumber\\
&=& \sum_{j=0}^{k}{k-j+\mb \choose k-j}L_{j}^{\nu}(x) 
- \sum_{j=0}^{k-1} {k-j+\mb-1 \choose k-j-1}L_{j}^{\nu}(x)+F_{k-1}(x).
\nonumber
\end{eqnarray}
Since $L_{0}^{\nu}(x)=1$ and $F_{0}(x)=L_{0}^{2\ma}(x)=1$,
\begin{eqnarray}
&&F_{k}(x)- \sum_{j=0}^{k}{k-j+\mb \choose k-j}L_{j}^{\nu}(x)
=F_{k-1}(x)-\sum_{j=0}^{k-1} {k-1-j+\mb \choose k-1-j}L_{j}^{\nu}(x)
\nonumber\\
&& \quad
= F_{0}(x)-{\mb \choose 0}L_{0}^{\nu}(x)=0.
\nonumber
\end{eqnarray}
Then we have (\ref{eqn:lemB2_1}).
From (\ref{eqn:F1}) and (\ref{eqn:G1}) we have
\begin{equation}
G_{k}(x)=-F_{k}(x)+\frac{k+2\ma}{k} F_{k-2}(x),
\quad k \in \N.
\label{eqn:G2}
\end{equation}
Then (\ref{eqn:lemB2_1}) gives (\ref{eqn:lemB2_2}).
This completes the proof.
\qed

\subsection{Proof of Lemma \ref{lem:skew0}}

We introduce a symmetric inner product
\begin{equation*}
(f, g) \equiv \int_{0}^{\infty} dx \, e^{-x} x^{2\ma} f(x) g(x).
\end{equation*}
It is easy to see that 
it is related with the elementary skew-symmetric
inner product (\ref{eqn:eskew1}) as
\begin{equation}
\label{eqn:skew=sym}
\langle f, g \rangle_{*}= 2({\cal I} f, g )
- \int_{0}^{\infty} dx \, e^{-x/2} x^{\ma}f(x)
\int_{0}^{\infty} dy \, e^{-y/2} y^{\ma}g(y),
\end{equation}
where
$$
{\cal I} f(x) \equiv \frac{\int_0^x dz \ e^{-z/2}z^\ma f(z)}
{e^{-x/2}x^\ma}.
$$
We consider the polynomials
\begin{eqnarray}
\label{def:Wn}
W_j(x) &=& L_{j}^{2\ma}(x)-\frac{j+2\ma}{j} L_{j-1}^{2\ma}(x),
\quad j\in\N.
\end{eqnarray}

\begin{lem}
\label{thm:lemB3}
For $j=\N_0, \; k\in\N$
\begin{eqnarray}
&&(W_k, F_j) = - \frac{\Gamma(j+2\ma+2)}{(j+1)!}\delta_{k-1 \ j}.
\label{eqn:WF}
\end{eqnarray}
\end{lem}

\noindent{\it Proof.} 
For $2\ma >-1$ the above equation can be derived 
immediately from 
the definitions (\ref{def:Wn}), (\ref{eqn:F2}),
and the orthogonal relation (\ref{eqn:orth}).
We can extend it to the case $\ma >-1$.
For $k=2,3,\dots$ it is easy to see that
\begin{eqnarray}
(L_k^{2\ma}, x^j) 
&=& (-1)^k \Gamma (k+2\ma+1)\delta_{j k},
\quad j=1,2,\dots,k.
\label{eqn:(L_n,x^n)}
\end{eqnarray}
However, if $2\ma< -1$, $(L_k^{2\ma}, 1)=\infty$.
Then we use the following equation:
\begin{equation}
(L_k^{2\ma}, L_j^{2\ma}-L_j^{2\ma}(0))
=\frac{\Gamma(k+2\ma+1)}{k!}\delta_{j k},
\quad k \in\N_0,
\label{(L,L-L0)}
\end{equation}
which is obtained from (\ref{eqn:(L_n,x^n)}) for $j=1,2,\dots, k$, 
since $L_k^{2\ma}(x)$ is a polynomial
whose coefficient of the $k$-th order is $(-1)^k/k!$.
By simple calculation we have
\begin{equation*}
(1, L_j^{2\ma}-L_j^{2\ma}(0))=0,
\quad j\in\N_0,
\label{(1,L-L0)}
\end{equation*}
and thus
\begin{equation*}
(L_k^{2\ma}, L_j^{2\ma}-L_j^{2\ma}(0))
=(L_k^{2\ma}-L_k^{2\ma}(0), L_j^{2\ma}-L_j^{2\ma}(0))=0,
\quad j \geq k+4.
\end{equation*}

From the definitions of $W_k$ and $F_0\equiv 1$ 
it is easily to see that
\begin{eqnarray}
(W_k, F_0) 
&=& -\Gamma(2\ma+2)\delta_{k-1 \ 0}.
\label{eqn:(W_1,1)}
\end{eqnarray}
When $k,j \in \N$,
from (\ref{eqn:(W_1,1)}) and (\ref{(L,L-L0)}),
\begin{eqnarray}
(W_k, F_j) 
&=& \left(
 L_k^{2\ma}-\frac{k+2\ma}{k}L_{k-1}^{2\ma}, \sum_{p=0}^j L_p^{2\ma}
\right)
\nonumber\\
&=& \sum_{p=1}^j \left(  
L_k^{2\ma}-\frac{k+2\ma}{k}L_{k-1}^{2\ma}, L_p^{2\ma}-L_p^{2\ma}(0)
\right)
\nonumber\\
&=& - \frac{\Gamma(k+2\ma+1)}{k!}\delta_{k-1 \ j}.
\end{eqnarray}
This completes the proof.
\qed

\vskip 3mm

Here we prove the following integral formula.

\begin{lem}
\label{thm:lemB1}
For $j\in\N, \ell\in\N_0$
\begin{eqnarray}
&& \int_{0}^{z} dx \, e^{-x/2} x^{\ma} G_j(x)
= 2 e^{-z/2} z^{\ma} W_j(x),
\label{eqn:G3}
\\
&& \int_{0}^{z} dx \, e^{-x/2} x^{\ma} F_{2\ell}(x) 
= 2^{\ma+1} {\ell +\ma \choose \ell}\gamma(\ma+1,z/2)
\nonumber\\
&& \qquad - 2 e^{-z/2} z^{\ma}{\ell +\ma \choose \ell}
\sum_{r=0}^{\ell-1} {\ell-r+\ma \choose \ell-r}^{-1}W_{2\ell-2r}(z),
\label{eqn:F3}
\end{eqnarray}
where $\gamma(c, y), c>0$ is the incomplete gamma function
$
\gamma(c, y) = \int_{0}^{y} dx \, e^{-x} x^{c-1}.
$

\end{lem}

\vskip 0.3cm
\noindent{\bf Remark 7. } \quad
If we set $\ma=0$ in (\ref{eqn:G3}), we will have 
the simpler equation
\begin{equation}
\int_{0}^{z} dx \, e^{-x/2} \frac{d}{dx} \Bigg\{
L_{j+1}^{0}(x)-L_{j-1}^{0}(x)
\Bigg\}
= 2 e^{-z/2} \Bigg\{ L_{j}^{0}(z) - L_{j-1}^{0}(z) \Bigg\},
\quad j\in\N.
\nonumber
\end{equation}


\noindent{\it Proof of Lemma \ref{thm:lemB1}.} \quad
We first introduce the functions defined by
\begin{eqnarray*}
\label{eqn:not1}
\psi_{j}^{2\ma}(z) &=& \int_{0}^{z} dx \,
e^{-x/2} x^{\ma} \frac{d}{dx} L_{j}^{2\ma}(x), 
\quad j\in\N_0,
\\
\label{eqn:not2}
\varphi_{j}^{2\ma}(z) &=& \int_{0}^{z} dx \,
e^{-x/2} x^{\ma} L_{j}^{2\ma}(x),
\quad j\in\N_0.
\end{eqnarray*}
Then (\ref{eqn:form2}) gives
\begin{equation}
\varphi_{j}^{2\ma}(z)= - \psi_{j+1}^{2\ma}(z)+\psi_{j}^{2\ma}(z).
\label{eqn:form2b}
\end{equation}
On the other hand, by (\ref{eqn:form1}),
\begin{eqnarray}
\psi_{j}^{2\ma}(z) 
&=&j \int_{0}^{z} dx \, e^{-x/2} x^{\ma-1} W_{j}(x).
\nonumber
\end{eqnarray}
Noting the assumption $\ma>-1$ and 
the fact that $W_j(x)={\cal O}(x)$, in $x\to 0$,
\begin{eqnarray}
&&2\ma \psi^{2\ma}_{j}(z) 
= 2j \left\{ e^{-z/2} z^\ma W_j(z)+\int_{0}^{z} dx \, x^{\ma}
\left( \frac{1}{2} e^{-x/2} W_j(x)
- e^{-x/2} \frac{d}{dx} W_j(x) \right)\right\}
\nonumber\\
&&\quad=
2j e^{-z/2} z^{\ma} W_j(z) 
+ j \varphi_{j}^{2\ma}(z)- (j+2\ma) \varphi_{j-1}^{2\ma}(z)
-2j \psi_{j}^{2\ma}(z) + 2(j+2\ma) \psi_{j-1}^{2\ma}(z).
\nonumber
\end{eqnarray}
Then we use (\ref{eqn:form2b}) to eliminate 
$\varphi_{j}^{2\ma}(z)$ and $\varphi_{j-1}^{2\ma}(z)$,
and the equality (\ref{eqn:G3}) is obtained from the relation
$$
\int_0^z dx \ e^{-x/2}x^\ma G_j(x) 
= \psi_{j+1}^{2\ma}(z)-\frac{j+2\ma}{j}\psi_{j-1}^{2\ma}(z).
$$
Note that (\ref{eqn:G3}) gives a recurrence relation for 
$\psi_{j}^{2\ma}(z)$,
$$
\psi_{j+1}^{2\ma}(z)= \frac{j+2\ma}{j}\psi_{j-1}^{2\ma}(z)
+ 2e^{-z/2}z^{\ma}W_j(z).
$$
It is solved as
\begin{eqnarray}
\psi_{2\ell+1}^{2\ma}(z)
&=& {\ell+\ma \choose \ell} \psi_{1}^{2\ma}(z)
+ 2 e^{-z/2} z^{\ma}{\ell+\ma \choose \ell}
\sum_{r=0}^{\ell-1} {\ell-r+\ma \choose \ell-r}^{-1}W_{2\ell-2r}(z).
\nonumber
\end{eqnarray}
Since $L_{1}^{2\ma}(x)=(1+2\ma)-x$,
\begin{eqnarray}
\psi_{1}^{2\ma}(z) &=& \int_{0}^{z}dx \, e^{-x/2} x^{\ma}
\frac{d}{dx} L_{1}^{2\ma}(z) 
= - \int_{0}^{z} dx \, e^{-x/2} x^{\ma} 
= -2^{\ma+1} \gamma(\ma+1, z/2).
\nonumber
\end{eqnarray}
Then we have (\ref{eqn:F3}).
\qed

\vskip 3mm

Then we can prove Lemma \ref{lem:skew0}.
From Lemma \ref{thm:lemB1} and two relations 
(\ref{eqn:G1}) and (\ref{eqn:skew=sym}) we have 
\begin{eqnarray*}
\langle F_{2q}, G_{2\ell+1} \rangle_* 
&=& -\langle G_{2\ell+1}, F_{2q} \rangle_* 
=-4 (W_{2\ell+1}, F_{2q}),
\\
\langle G_{2q+1}, G_{2\ell+1} \rangle_* 
&=& 4 (W_{2q+1}, G_{2\ell+1})
\\
&=& 4 \left\{-(W_{2q+1}, F_{2\ell+1})
+\frac{2\ell+1+2\ma}{2\ell+1}(W_{2q+1}, F_{2\ell-1})\right\}.
\end{eqnarray*}
Then (\ref{eqn:sskew2}) and (\ref{eqn:44_2}) 
are derived from Lemma \ref{thm:lemB3}.
From Lemma \ref{thm:lemB1} and (\ref{eqn:WF}),  
\begin{eqnarray*}
&&\int_{0}^{\infty} dw \, e^{-w/2} w^{\ma} \gamma(\ma+1, w/2)F_{2\ell}(w)
\\
&&= \left[ \gamma(\ma+1, w/2)\int_{0}^w dz \, e^{-z/2} z^{\ma} F_{2\ell}(z)
\right]_0^\infty
\\
&&\quad - \int_0^\infty dw \ 2^{-(1+\ma)}e^{-w/2}w^\ma 
\int_{0}^w dz \, e^{-z/2} z^{\ma} F_{2\ell}(z)
\\
&&= 2^{\ma+1}{\ell+\ma \choose \ell}
\left\{ \Gamma (\ma+1)^2
- \int_0^\infty dw \ 2^{-(1+\ma)}e^{-w/2}w^\ma \gamma(\ma+1, w/2)
\right\}
\\
&&= 2^{\ma+1}{\ell+\ma \choose \ell}
\left\{ \Gamma (\ma+1)^2
- \frac{1}{2} \left[ \gamma(\ma+1, w/2)^2 \right]_0^\infty 
\right\}
= 2^{\ma}{\ell+\ma \choose \ell} \Gamma (\ma+1)^2.
\end{eqnarray*}
Then by using Lemma \ref{thm:lemB1} and (\ref{eqn:WF}) we have
(\ref{eqn:44_1}), since
\begin{eqnarray*}
\langle F_{2q}, F_{2\ell} \rangle_{*}
&=& 2 \int_{0}^{\infty} dw \, e^{-w/2} w^{\ma} F_{2\ell}(w)
\int_{0}^{w}dz \, e^{-z/2} z^{\ma} F_{2q}(z) 
\\
&&\qquad 
- \int_{0}^{\infty} dw \, e^{-w/2} w^{\ma} F_{2\ell}(w)
\int_{0}^{\infty}dz \, e^{-z/2} z^{\ma} F_{2q}(z) 
\\
&=& 2^{\ma+2}{q+\ma \choose q}
\int_{0}^{\infty} dw \, e^{-w/2} w^{\ma} \gamma(\ma+1, w/2)F_{2\ell}(w)
\\
&&\qquad -2^{2\ma+2} {q+\ma \choose q}{\ell+\ma \choose \ell}\Gamma (\ma+1)^2
= 0.
\end{eqnarray*}
This completes the proof.
\qed

\SSC{Inverse of $\{\alpha_{k,j}\}$}\label{chap:beta}

We put
\begin{equation}
Q_{2\ell}(x)=F_{2\ell}(x) 
\qquad\mbox{and}\qquad 
Q_{2\ell+1}(x)=G_{2\ell+1}(x),
\label{def:H}
\end{equation}
for $\ell \in\N_0$,
and  $Q_k\equiv 0$ for $k\in \Z_{-}$.
Lemma \ref{thm:lemB2} gives the expansion formula
of $Q_{k}(x)$ in terms of $\{ L_{j}^{\nu}(x) \}$.
Here we give the formula to expand
$L_{j}^{\nu}(x)$ in terms of $\{Q_{k}(x)\}$.
In other words, we calculate the inverse of
the matrix $\valpha =(\alpha_{k,j})$  given by
(\ref{def:alpha}), which is denoted by
$\vbeta = (\beta_{j, k})$ and used in Section 5.2.

Let
$
b(n)=(n+2\ma)/n, n\in\N,
$
and 
\begin{equation}
b(m,n)=
\left\{
   \begin{array}{ll}
      b(m)b(m+2)\cdots b(n),
 & \mbox{if $m,n$ are odd and $m\le n$}, \\
      1,
 & \mbox{if $m,n$ are odd and $m>n$}, \\
       0,
 & \mbox{otherwise}. \\
   \end{array}\right.
\label{eqn:bmn}
\end{equation}
Then the following lemma holds.

\vskip 3mm
\begin{lem}
\label{thm:lemC2}
For $j\in\N_0$
\begin{equation}
L_j^{\nu}(x) = \sum_{k=0}^j \beta_{j,k}Q_k(x).
\label{eqn:eq2}
\end{equation}
where $\beta_{j,k}$, $k,1,\dots,j$, $j\in\N_0$ are 
defined by the following:

\noindent When $k$ is even
\begin{equation}
\beta_{j,k}=
\left\{
   \begin{array}{ll}
      0,
 & \mbox{if} \ j<k, \\
 \displaystyle{{j-k-\mb-2 \choose j-k}},
 & \mbox{if} \ j\ge k, \\
   \end{array}\right.
\label{def:beta_even}
\end{equation}
and, when $k$ is odd
\begin{equation}
\beta_{j,k}=
\left\{
   \begin{array}{ll}
      0,
 & \mbox{if} \ j<k, \\
 - \displaystyle{\sum_{r=[(k+1)/2]}^{[(j+1)/2]}
 b(k+2,2r-1){j-2r-\mb-1 \choose j-2r+1}},
 & \mbox{if} \ j\ge k. \\
   \end{array}\right.
\label{def:beta_odd}
\end{equation}
\end{lem}

\noindent{\it Proof.}
From (\ref{eqn:F2}) and (\ref{eqn:G2})  we have
\begin{eqnarray*}
&&Q_{2\ell}(x)= \sum_{j=0}^{2\ell} L_j^{2\ma}(x),
\quad
Q_{2\ell+1}(x)= -\sum_{j=0}^{2\ell+1} L_j^{2\ma}(x)
+b(2\ell+1)\sum_{j=0}^{2\ell-1} L_j^{2\ma}(x).
\end{eqnarray*}
By simple calculations we see that
\begin{equation}
L_j^{2\ma}(x) = \sum_{k=0}^j \overline{\beta}_{j,k}Q_k(x),
\quad j\in\N_0
\label{eqn:L=bH}
\end{equation}
where 
$\overline{\beta}_{j,k}$, $k=0,1,\dots,j$, $j\in\N_0$,
are defined by the following:

\noindent When $k$ is even
\begin{equation*}
\overline{\beta}_{j,k}=
\left\{
   \begin{array}{ll}
      1,
 & \mbox{if} \ j=k, \\
      -1
 & \mbox{if} \ j=k+1, \\
 	  0,
 & \mbox{otherwise}, \\
   \end{array}\right.
\label{eqn:relation_even}
\end{equation*}
\noindent
and, when $k$ is odd
\begin{equation*}
\overline{\beta}_{j,k}=
\left\{
   \begin{array}{ll}
      b(k+2, j-1),
 & \mbox{if $j$ is even}, \\
      -b(k+2,j),
 & \mbox{if $j$ is odd}. \\
   \end{array}\right.
\label{eqn:relation_odd}
\end{equation*}
Using the formula (\ref{eqn:form3}),
(\ref{eqn:L=bH}) gives
\begin{eqnarray}
L_j^{\nu}(x) 
&=& \sum_{p=0}^j {j-p-\mb-1 \choose j-p}
\sum_{k=0}^p \overline{\beta}_{p,k}Q_k(x).
\nonumber
\end{eqnarray}
If we define 
\begin{equation}
\beta_{j,k}= \sum_{p=k}^j \overline{\beta}_{p,k}{j-p-\mb-1 \choose j-p}
\label{eqn:AppCeq}
\end{equation}
for $j \ge k$ and $\beta_{j,k}=0$ for $j < k$,
(\ref{eqn:eq2}) is satisfied.
The expressions 
(\ref{def:beta_even}) and
(\ref{def:beta_odd}) are derived from
(\ref{eqn:AppCeq}) by simple calculation.
\qed

\SSC{Elementary Calculation for Asymptotics Estimation}
\label{chap:Elementary}

By Stirling's formula
$\Gamma(x) \sim \sqrt{2\pi}x^{x-1/2}e^{-x},
\, x\to\infty$,
we have 
\begin{eqnarray}
\frac{\Gamma(n+\alpha+1)}{\Gamma(n+1)}
&\sim& (n+1)^\alpha,
\quad n\to\infty,
\label{sim:Gamma/Gamma}
\\
{n+\alpha \choose n}
&\sim& \frac{(n+1)^\alpha}{\Gamma(\alpha+1)}
\quad n\to\infty,
\label{sim:choose}
\end{eqnarray}
for any $\alpha \in \R\setminus \Z_-$, and 
\begin{eqnarray}
&&b(1, 2\ell-1)=
\frac{\Gamma (2\ma+1)}{2^\ma \Gamma (\ma+1)}
{2\ell+2\ma \choose 2\ell}/{\ell+\ma \choose \ell}
\sim (2\ell)^\ma, \quad \ell\to\infty
\label{sim:b1l}
\\
&&b(2\ell+1, 2p-1)=
\frac{b(1,2p-1)}{b(1,2\ell-1)}\sim  \left(\frac{p}{\ell}\right)^{\ma},
\quad \ell\to\infty
\label{sim:blp}
\end{eqnarray}
for $\ell, p\in\N$ with $\ell < p$,
where $b(m,n)$ is defined by (\ref{eqn:bmn}).

From now on, we assume that $T=N$, $t_m=T+s_m$ with
$s_{m} < 0$.
We set
\begin{equation}
2 \ell = N \theta \quad \mbox{and} \quad 
j=2 \ell \eta,
\label{eqn:setting}
\end{equation}
and consider the limit $N \to \infty$
with some $ \eta, \theta \in (0, \infty)$.
Then we have
\begin{equation*}
\chi_m^j= \left( \frac{2T-t_m}{t_m} \right)^j
=\left( 1-\frac{2s_m}{t_m} \right)^{N\theta\eta}
\sim \exp \left( -2s_m \theta \eta \right),
\quad N\to\infty,
\label{sim:chi}
\end{equation*}
and
\begin{equation}
\sum_{p=0}^\alpha (-1)^p {\alpha \choose p}\chi_m^{j-p}
\sim (2\ell)^{-\alpha} 
\left( \frac{d}{d\eta}\right)^{\alpha}
\exp \left( -2s_m \theta \eta \right),
\quad N\to\infty.
\label{sim:chi_2}
\end{equation}

We use the following identities
(see (\ref{def:combination}) and
pages 8, 201 and 202 in \cite{Rio79}).


\noindent {\rm (i)} \quad
Let $\alpha \in\N_0$ and $c\in\R$. Then
\begin{equation}
\sum_{p=0}^\alpha (-1)^p {\alpha \choose p}{n-p+c \choose n-p-j}
={n-\alpha+c \choose n-j}.
\label{formula:combi_1}
\end{equation}

\noindent {\rm (ii)} \quad
Let $\alpha \in\N_0$, $c\in\R$
and $a_k$, $k=1,2,\dots,$ be a sequence in $\R$.
Then
\begin{eqnarray}
&&\sum_{r=0}^{\infty}{r+c \choose r}a_r
= \sum_{r=0}^{\infty}{r+c+\alpha \choose r}
\sum_{p=0}^{\alpha} (-1)^p {\alpha \choose p} a_{r+p},
\label{formula:combi_2}
\end{eqnarray}

\noindent {\rm (iii)} \quad
Let $\alpha \in \N_0$,
and $a_k, b_k$, $k=1,2,\dots,$ be sequences in $\R$.
Then
\begin{eqnarray}
&&\sum_{k=0}^{\alpha}(-1)^k {\alpha \choose k} a_k b_k
= \sum_{\beta=0}^{\alpha}{\alpha \choose \beta}
\sum_{p=0}^{\beta} (-1)^p {\beta \choose p} a_p 
\sum_{q=0}^{\alpha -\beta}(-1)^q {\alpha-\beta \choose q} b_q.
\label{formula:combi_3}
\end{eqnarray}

\vskip 3mm 

\begin{lem}
\label{thm:lemD1}
\quad For any $\alpha\in\N_0$ and $w\ge 0$ we have
\begin{equation*}
\sum_{p=0}^\alpha (-1)^p {\alpha \choose p}
L_{n-p}^{\nu}\left(\frac{w}{n}\right)
\sim
\left(\frac{w}{n}\right)^{\alpha-\nu}
\left(\frac{d}{dw}\right)^\alpha \left\{ w^{\nu/2} J_{\nu}(2\sqrt{w})\right\},
 \quad n\to\infty.
\label{eqn:L-L=dJ}
\end{equation*}
\end{lem}

\noindent{\it Proof. } 
From the definition of the
Laguerre polynomials (\ref{eqn:Laguerre}) 
and (\ref{formula:combi_1}), we have
\begin{eqnarray*}
\sum_{p=0}^\alpha (-1)^p {\alpha \choose p}L_{n-p}^{\nu}\left(y \right)
&=& \sum_{p=0}^\alpha (-1)^p {\alpha \choose p}
\sum_{j=0}^{n-p} (-1)^j {n-p+\nu \choose n-p-j}\frac{y^j}{j!}
\nonumber\\
&=& \sum_{j=0}^{n} \frac{(-1)^j}{j!}{n-\alpha+\nu \choose n-j}y^j.
\end{eqnarray*}
Hence, by (\ref{sim:choose})
\begin{eqnarray*}
\lim_{n\to\infty}
\left(\frac{w}{n}\right)^{\nu-\alpha}
\sum_{p=0}^\alpha (-1)^p 
{\alpha \choose p}L_{n-p}^{\nu}\left(\frac{w}{n}\right)
&=& \lim_{n\to\infty}
\sum_{j=0}^{n} \frac{(-1)^j}{j!}{n-\alpha+\nu \choose n-j}
\left(\frac{w}{n}\right)^{j+\nu-\alpha}
\nonumber\\
&=&\sum_{j=0}^{\infty}(-1)^j
\frac{w^{\nu+j-\alpha}}{\Gamma(j-\alpha+\nu+1)j!}
\nonumber\\
&=&\left( \frac{d}{dw}\right)^\alpha
\left( \sum_{j=0}^{\infty}(-1)^j
\frac{w^{\nu+j}}{\Gamma(\nu+j+1) j!}\right)
\end{eqnarray*}
Then we obtain the lemma
\qed
\vskip 3mm

Applying the above lemma, we obtain 
the following asymptotics,
where $L_{j}^{\nu}$ and $\hL_{j}^{\nu}$
are defined by (\ref{eqn:LLhat}).

\begin{lem}
\label{thm:lemD2}
\quad For any $\alpha\in\N_0$, 
$\theta, \eta \in (0, \infty)$, and $x \in \Rp$
we have 
\begin{eqnarray}
&&\sum_{p=0}^\alpha 
(-1)^p {\alpha \choose p}
L_{j-p}^{\nu}\left(\frac{x}{N}, -s_m \right)
\sim \frac{(2\ell)^{\nu-\alpha}}{(\theta x)^\nu}
\Knu^{(\alpha)}(\theta,\eta,x,-s_m), \quad N \to \infty,
\label{sim:K}
\\
&&\sum_{p=0}^\alpha (-1)^p {\alpha \choose p}
\hL^{\nu}_{j+p}\left(\frac{x}{N}, s_{m} \right)
\sim (2\ell)^{-\alpha}
\hKnu^{(\alpha)}(\theta,\eta,x,s_m), \quad N \to \infty.
\label{sim:hK}
\end{eqnarray}
\end{lem}

\noindent{\it Proof. }
Since 
\begin{eqnarray*}
&&\sum_{p=0}^\alpha (-1)^p {\alpha \choose p}
L^{\nu}_{j-p}\left(\frac{x}{N}\right)\chi_m^{j-p}
\\
&&=\sum_{\beta=0}^\alpha {\alpha \choose \beta}
\sum_{p=0}^\beta (-1)^p {\beta \choose p}
L^{\nu}_{j-p}\left(\frac{x}{N}\right)
\sum_{q=0}^{\alpha-\beta} (-1)^q {\alpha-\beta \choose q}
\chi_m^{j-q},
\end{eqnarray*}
by (\ref{formula:combi_3}), the asymptotic (\ref{sim:K})
is derived from (\ref{sim:chi_2}) and Lemma \ref{thm:lemD1}
with $n=j=N \theta \eta, w=\theta\eta x$.
From (\ref{formula:combi_3}), we have 
\begin{eqnarray}
&&\sum_{k=0}^\alpha (-1)^k {\alpha \choose k}
\hL^{\nu}_{j+k}\left(\frac{x}{N}, s_m \right)
\nonumber\\
&&= \sum_{\beta=0}^{\alpha}{\alpha \choose \beta}
\sum_{p=0}^{\beta} (-1)^p {\beta \choose p} 
L^{\nu}_{j+p}\left( \frac{x}{N} \right)\chi_m^{-(j+p)}
\sum_{q=0}^{\alpha -\beta}(-1)^q {\alpha-\beta \choose q}
\frac{\Gamma(j+q+1)}{\Gamma(j+q+1+\nu)}.
\nonumber\\
\label{eqn:L-Gamma}
\end{eqnarray}
By (\ref{sim:Gamma/Gamma}), we see
\begin{eqnarray*}
\sum_{q=0}^{\alpha-\beta}(-1)^q {\alpha-\beta \choose q} 
\frac{\Gamma (j+q+1)}{\Gamma (j+q+1+\nu)}
&=& \nu \sum_{q=0}^{\alpha-\beta-1}(-1)^q {\alpha-\beta-1 \choose q} 
\frac{\Gamma (j+q+1)}{\Gamma (j+q+2+\nu)}
\\
&=& 
\frac{\nu(\nu+1)\cdots(\nu+\alpha-\beta-1)\Gamma (j+1)}
{\Gamma (j+1+\alpha-\beta+\nu)},
\\
&\sim& (2\ell)^{-(\alpha-\beta+\nu)}
\left(-\frac{d}{d\eta}\right)^{\alpha-\beta}\eta^{-\nu}, 
\quad N \to \infty.
\end{eqnarray*}
On the other hand, (\ref{sim:K}) gives
\begin{eqnarray*}
&&\sum_{p=0}^{\beta} (-1)^p {\beta \choose p} 
L^{\nu}_{j+p}\left(\frac{x}{N} \right)\chi_m^{-(j+p)}
\sim (2\ell)^{\nu-\beta}\left(- \frac{d}{d\eta}\right)^{\beta}
\left\{ \eta^{\nu} \hKnu (\theta,\eta,x,s_m)\right\}.
\end{eqnarray*}
Hence, the asymptotic (\ref{sim:hK})
is derived from (\ref{eqn:L-Gamma}).
\qed

\SSC{On Temporally Homogeneous Limit}\label{chap:AppLimit}

\begin{lem}
\label{thm:lemE1}
\quad For any $c\in \R$ and $\eta, \theta,x \ge 0$, 
we have that as $t\to\infty$
\begin{eqnarray}
&&\Knu^{(c)}(\theta,\eta,x,t) 
\sim (2t\theta)^c(\theta\eta x)^{\nu/2}
J_\nu (2\sqrt{\theta\eta x})e^{2t\theta\eta},
\label{sim:k^c(t)}
\\
&&\hKnu^{(c)}(\theta,\eta,x,-t) 
\sim (2t\theta)^c(\theta\eta x)^{-\nu/2}
J_\nu (2\sqrt{\theta\eta x})e^{-2t\theta\eta}, 
\label{sim:hk^c(t)}
\\
&&\int_1^\infty d\xi \ \xi^\ma \hKnu^{(c+1)}(\theta,\xi,x,-t)
\sim (2t\theta)^{c} 
(\theta x)^{-\nu/2}
J_\nu (2\sqrt{\theta x})e^{-2t\theta}.
\label{sim:intk^c(t)}
\end{eqnarray}
\end{lem}

\noindent{\it Proof. }
From the expression (\ref{def:K^c})
with the definition (\ref{def:k_nu}), we have
\begin{eqnarray*}
&&\Knu^{(c)}(\theta,\eta,x,t) 
=\frac{e^{2t\theta\eta}}{\Gamma(-c)}
\sum_{k=0}^\infty \frac{(-1)^k \eta^{k-c}}{k! (k-c)}
\sum_{j=0}^k {k \choose j}
\Knu^{(j)}(\theta,\eta,x,0)
(2t\theta)^{k-j}
\\
&&\qquad\qquad\qquad
\sim \frac{e^{2t\theta\eta}}{\Gamma(-c)}
(2t\theta)^{c}\Knu (\theta,\eta,x,0)
\sum_{k=0}^\infty
\frac{(-1)^k (2t\theta\eta)^{k-c}}{k! (k-c)},
\quad t\to\infty.
\end{eqnarray*}
From the relation
$$
\frac{d}{dz} \sum_{k=0}^\infty
\frac{(-1)^k z^{k-c}}{k! (k-c)}
= z^{-c-1}e^{-z},
$$
and the equation 
$$
\Gamma(-c)= \sum_{k=0}^\infty\frac{(-1)^k}{k! (k-c)}
+\int_1^\infty dz \ z^{-c-1}e^{-z}
$$
(see (1.1.19) in \cite{AAR99}),
we have
$$
\Gamma(-c)= \lim_{t\to\infty}\sum_{k=0}^\infty
\frac{(-1)^k (2t\theta\eta)^{k-c}}{k! (k-c)}.
$$
Then we conclude 
\begin{eqnarray*}
&&\Knu^{(c)}(\theta,\eta,x,t) 
\sim e^{2t\theta\eta}(2t\theta)^{c}\Knu (\theta,\eta,x,0)
=(2t\theta)^{c}\Knu (\theta,\eta,x,t),
\quad t \to \infty.
\end{eqnarray*}
Hence (\ref{sim:k^c(t)}) is derived from (\ref{def:k_nu}).

Let $n=[c+1]_{+}$ and $\beta>0$ with $c=n-\beta$.
Since
$$
\left(- \frac{d}{d\eta}\right)^n \hKnu(\theta,\eta,x,-t)
\sim (2t\theta)^{n}\hKnu(\theta,\eta,x,-t),
\quad t\to\infty,
$$
(\ref{def:hK^c}) gives
\begin{eqnarray*}
&&\hKnu^{(c)}(\theta,\eta,x,-t)
\sim \frac{(2t\theta)^{n}}{\Gamma (\beta)}
\int_0^\infty d\xi \ \xi^{\beta-1} \hKnu (\theta,\eta+\xi,x,-t)
\\
&&\qquad\qquad\qquad
=\frac{(2t\theta)^{n-\beta}}{\Gamma (\beta)}
\int_0^\infty d\zeta \ \zeta^{\beta-1} 
\hKnu (\theta,\eta+\frac{\zeta}{2t\theta},x,-t)
\\
&&\qquad\qquad\qquad
\sim (2t\theta)^{c}
\hKnu (\theta,\eta,x,-t),
\quad t \to \infty.
\end{eqnarray*}
Then (\ref{sim:hk^c(t)}) is derived from (\ref{def:hk_nu}).

From (\ref{sim:hk^c(t)}) we have
\begin{eqnarray*}
\int_1^\infty d\xi \ \xi^\ma \hKnu^{(c+1)}(\theta,\xi,x,-t)
&\sim& (2t\theta)^{c+1}
\int_1^\infty d\xi \ \xi^\ma \hKnu (\theta,\xi,x,0)e^{-2t\theta\xi}
\\
&\sim& (2t\theta)^{c} \hKnu (\theta,1,x,-t),
\quad t \to \infty.
\end{eqnarray*}
This completes the proof.
\qed
\vskip 3mm

Applying the above lemma, we have as $s_m, s_n \to -\infty$
with the difference $s_n - s_m$ fixed
\begin{eqnarray*}
\cD (s_m,x; s_n,y)
&\sim& \frac{(xy)^{\mb/2}(s_n-s_m)}{2^{2\mb+3}(s_ms_n)^{\mb+1}}
\int_0^1 d\theta \ \theta^{-\mb}
J_{\nu}(2\sqrt{\theta x})J_{\nu}(2\sqrt{\theta y})
e^{-2(s_m+s_n)\theta}
\nonumber\\
&\sim& 
\frac{(xy)^{\mb/2}(s_m-s_n)}{2^{2\mb+4}(s_m+s_n)(s_ms_n)^{\mb+1}}
J_{\nu}(2\sqrt{x})J_{\nu}(2\sqrt{y})e^{-2(s_m+s_n)},
\end{eqnarray*}
\begin{eqnarray*}
\cI(s_m,x; s_n,y)
&\sim& 
\frac{2^{2\mb+1}(s_ms_n)^{\mb}(s_n-s_m)}{(xy)^{\mb/2}}
\int_1^\infty d\theta \ \theta^{\mb}
J_{\nu}(2\sqrt{\theta x})J_{\nu}(2\sqrt{\theta y})
e^{2(s_m+s_n)\theta}
\nonumber\\
&\sim& 
\frac{2^{2\mb}(s_ms_n)^{\mb}(s_n-s_m)}{(s_m+s_n)(xy)^{\mb/2}}
J_{\nu}(2\sqrt{x})J_{\nu}(2\sqrt{y})e^{2(s_m+s_n)},
\end{eqnarray*}
and
\begin{eqnarray*}
&&\cS (s_m, x; s_n,y)
\sim 
\left(\frac{y}{x}\right)^{\mb/2}
\int_0^1 d\theta \ 
J_{\nu}(2\sqrt{\theta x})J_{\nu}(2\sqrt{\theta y})
e^{2(s_m-s_n)\theta}.
\end{eqnarray*}
It is then clear that
\begin{equation*}
\lim_{s_m, s_n\to -\infty} \cD(s_m,x;s_n,y)\cI(s_m,x;s_n,y)=0.
\end{equation*}

\section*{Acknowledgment} 

The authors would like to thank T. Sasamoto 
and T. Imamura 
for useful discussion on determinantal and
Pfaffian processes.


\end{document}